\documentclass[12pt]{amsart}
\usepackage{amsthm}
\usepackage{bbm}
\usepackage{hyperref}
\newtheorem{theorem}{Theorem}[section]
\newtheorem{lemma}[theorem]{Lemma}
\newtheorem{proposition}[theorem]{Proposition}
\theoremstyle{definition}

\newtheorem{corollary}[theorem]{Corollary}
\theoremstyle{remark}
\newtheorem{remark}[theorem]{Remark}

\usepackage{tikz}
\usetikzlibrary{calc}

\counterwithin*{section}{part}

\usepackage{amssymb}
\usepackage{empheq}
\usepackage{stmaryrd}
\usepackage{enumerate}
\usepackage[shortlabels]{enumitem}
\usepackage{calc}
\usepackage{url}
\usepackage{comment}
\usepackage{mathabx}

\newcommand{\norm}[1]{\left\lVert#1\right\rVert}
\usepackage{mathtools}

\def\EE{\mathbb{E}}
\def\PP{\mathbb{P}}
\def\NN{\mathbb{N}}
\def\ZZ{\mathbb{Z}}
\def\QQ{\mathbb{Q}}

\newcommand{\saw}[1]{\mathsf{SAW}_{#1}}

\newcommand{\psaw}[1]{\mathbb{P}_{\saw{n}}}
\newcommand{\isab}{\mathsf{iSAB}}
\newcommand{\pisab}{\mathbb{P}_{\mathsf{iSAB}}}
\newcommand{\eisab}{\mathbb{E}_{\mathsf{iSAB}}}
\newcommand{\pisabn}{\mathbb{P}^{\otimes\NN}_{\mathsf{iSAB}}}
\newcommand{\eisabn}{\mathbb{E}^{\otimes\NN}_{\mathsf{iSAB}}}

\newcommand{\R}{\mathbb{R}}
\newcommand{\Z}{\mathbb{Z}}
\newcommand{\N}{\mathbb{N}}

\newcommand{\SAW}{\mathsf{SAW}}
\newcommand{\SAB}{\mathsf{SAB}}
\newcommand{\iSAB}{\mathsf{iSAB}}
\newcommand{\ZigZag}{\mathsf{ZZ}}
\newcommand{\shortZZ}{\mathsf{shortZZ}}
\newcommand{\manyZZfewR}{\mathsf{ManyZZfewR}}

\renewcommand{\P}{\mathbb{P}}
\newcommand{\irrP}{\P_{\iSAB}}
\newcommand{\infP}{\P_{\iSAB}^{\otimes \N}}

\newcommand{\E}{\mathbb{E}}
\newcommand{\irrE}{\E_{\iSAB}}
\newcommand{\infE}{\E_{\iSAB}^{\otimes \N}}
\DeclareMathOperator{\spn}{\mathbf{span}}

\usepackage[left=2.0cm,%
                right=2.0cm,%
                top=2.5cm,%
                bottom=3.5cm,%
                headheight=12pt,%
                a4paper]{geometry}%

\tikzset{
walk/.style={black,line width=0.85pt,line cap=round,line join=round},
cone/.style={gray!72,densely dashed,line width=0.55pt},
point/.style={circle,fill=blue,draw=white,line width=0.3pt,inner sep=1.65pt},
chosen/.style={circle,fill=red,draw=white,line width=0.3pt, inner sep=1.85pt},
past/.style={line width=1.1pt,line cap=round,line join=round},
future/.style={line width=1.1pt,line cap=round,line join=round}
}

\begin{document}

\title[Quantitative sub-ballisticity of self-avoiding walk]{Quantitative sub-ballisticity of self-avoiding walk on $\mathbb{Z}^d$}

\author[Ekin Arik\"ok, Christoforos Panagiotis, Nicholas Pischke, Thomas Powell]{Ekin Arik\"ok${}^{\MakeLowercase a}$, Christoforos Panagiotis${}^{\MakeLowercase b}$, Nicholas Pischke${}^{\MakeLowercase c}$, Thomas Powell${}^{\MakeLowercase c}$}
\date{\today}
\maketitle
\vspace*{-5mm}
\begin{center}
{\scriptsize 
${}^a$ Section de Math\'ematiques, Universit\'e de Gen\`eve,\\
Rue du Conseil-G\'en\'eral 7-9, 1205 Geneva, Switzerland,\\
${}^b$ Department of Mathematics, University of Bath,\\
Claverton Down, Bath, BA2 7AY, United Kingdom,\\
${}^c$ Department of Computer Science, University of Bath,\\
Claverton Down, Bath, BA2 7AY, United Kingdom,\\
E-mails: ekin.arikoek@unige.ch, cp2324@bath.ac.uk, nnp39@bath.ac.uk, trjp20@bath.ac.uk}
\end{center}

\begin{abstract}
We prove quantitative sub-ballisticity for the self-avoiding walk on the lattice $\mathbb{Z}^d$, $d\geq 2$. Namely, we obtain an explicit exponential rate of decay for the probability that a self-avoiding walk of length $n$ reaches distance at least $vn$, in terms of $v>0$. As a consequence, we show that the displacement of an $n$-step self-avoiding walk is at most $n/\exp(c\log^* n)$ with very high probability. Previously, only a non-quantitative $o(n)$ bound was known from the work of Duminil-Copin and Hammond. Along the way, we show that at criticality the partition function of bridges of length $n$ is on average bounded by $C/\log\log n$, which we believe to be of independent interest. Our approach is based on a quantitative analysis of the original nonconstructive proof of Duminil-Copin and Hammond using ideas from proof theory.
\end{abstract}

\section{Introduction}

A \emph{self-avoiding walk} on a graph is a walk that visits each vertex of the graph at most once. The model arose in the study of polymer chains \cite{Flory53}, and has since attracted the interest of mathematicians and physicists alike.  We refer to \cite{LecturesSAW,MaSaSAWBook} for comprehensive introductions. In this article, we are concerned with the self-avoiding walk on the hypercubic lattice $\mathbb{Z}^d$, and specifically on producing new quantitative results on their \emph{displacement}. Our main result is achieved through an analysis of a nonconstructive proof of the sub-ballisticity of self-avoiding walks due to Duminil-Copin and Hammond \cite{DCHammond2013}. Ideas from the logic-based \emph{applied proof theory} programme \cite{Kohlenbach2008} formed part of our methodology, and indeed, the broad structure of this analysis can be explained in terms of the Dialectica interpretation of G\"odel \cite{Goedel1958}. Though our main result can be understood without reference to mathematical logic, we consider this connection between self-avoiding walks and proof theory to be of independent interest.

\subsection{Background: Self-avoiding walks in $\ZZ^d$} 

In dimensions $d\geq 2$ the number of self-avoiding walks of length $n$ starting from $0$, which we denote by $c_n$, grows exponentially fast and has a well-defined exponential growth rate
\[
\mu_c:=\lim_{n\to\infty}c_n^{1/n}
\]
that we call the \emph{connective constant} of the graph. The exact value of $\mu_c$ is not known in any dimension $d\geq 2$ and to the best of our knowledge is not expected to take an interesting value. In contrast, on the hexagonal lattice, Duminil-Copin and Smirnov \cite{DC12} proved that $\mu_c=\sqrt{2+\sqrt{2}}$ by introducing and studying the so-called parafermionic observable.  
    
The main questions in the study of self-avoiding walk concern the asymptotic growth of the number of self-avoiding walks of length $n$ starting from a fixed vertex and the geometry of a typical long walk. An important direction in the early study of the model was to obtain quantitative improvements on the bound $c_n\leq e^{o(n)} \mu_c^n$. Hammersley and Welsh \cite{HammersleyWelsh} introduced the notion of \emph{bridges}, a special class of self-avoiding walks, and used an unfolding operation to compare $c_n$ with the number of bridges of length $n$ starting from $0$, which implies that there exists a constant $C>0$ such that for every $n\geq 1$,
\[
c_n\leq e^{C\sqrt{n}} \mu_c^n.
\]
See also \cite{Hammersley1961, KestenI, KestenII} for related works.

Despite this initial progress, improving on these early results has proved difficult. A more complete picture has, however, emerged in high dimensions. As expected for many models in statistical mechanics, in high-dimensions the self-avoiding walk exhibits mean-field behaviour, in the sense that the model forgets the geometry of the lattice and the critical exponents that describe its large-scale behaviour take the same values as those on a tree or on the complete graph. This was established by Hara and Slade \cite{HaSla1992,HaSlaCri}, who used the lace expansion to prove, in particular, that self-avoiding walk converges to Brownian motion under diffusive scaling in dimensions $d\geq 5$. In contrast, in low dimensions, the gap between what is known and what is expected to be true is very large. See \cite{PolyIns, SAWEnds, DCKY14, DCHammond2013, GM20, Hammond18,Madras14} and references therein for some relevant results.

\subsection{Main results} 

The focus of this paper is on the displacement of self-avoiding walks. Duminil-Copin and Hammond \cite{DCHammond2013} have shown that the self-avoiding walk in $\ZZ^d$ is sub-ballistic in the following sense. Fix $d\geq 2$ and let $\saw{n}$ denote the set of self-avoiding walks of length $n$ on $\ZZ^d$ that start at $0$. Let $\psaw{n}$ denote the uniform law on $\saw{n}$. 

\begin{theorem}
[Theorem 1.1 of \cite{DCHammond2013}]
\label{res:main}
For any $v>0$ there exists $\varepsilon>0$ such that, for every $n \in \NN$ large enough,
\[
\psaw{n}\left(\max\left\{\norm{\gamma_k}\mid 0\leq k\leq n\right\}\geq vn\right)\leq e^{-\varepsilon n}.
\]
\end{theorem}
In the above theorem, the value of $\varepsilon$ is ineffective, so the theorem only implies that a typical self-avoiding walk of length $n$ does not exit the ball of radius $R_n$ for some non-quantitative $R_n=o(n)$. The main result of this paper is an explicit construction of the exponential rate of decay $\varepsilon$ in terms of $v$, which in particular gives that a typical self-avoiding walk of length $n$ does not exit the ball of radius $R_n=O(n/\exp(c\log^*n))$.
\begin{theorem}[Main result]
\label{thm: main theorem}
For any $v>0$ and any $n \geq n_0(v)$,
\[
\psaw{n}\left(\max\left\{\norm{\gamma_k}\mid 0\leq k\leq n\right\}\geq vn\right)\leq e^{-n/T(v)}
\]
where $T(v) = \exp^{\circ \, O(\log(1/v))}(1)$ and $n_0(v) = 4C_{HW}^2T(v)^2+2T(v)\log(1+2d)$, with $C_{HW}>0$ the constant defined in Proposition \ref{HW}.
\end{theorem}

\begin{corollary}
    There exists $c > 0$ and $n_0 \in \N$ such that for each $n \geq n_0$,
\[
\psaw{n}\left(\max\left\{\norm{\gamma_k}\mid 0\leq k\leq n\right\}\geq \frac{n}{\exp(c\log^\ast n)}\right)\leq e^{-n^{2/3}}. 
\]
Moreover, $c>0$ and $n_0$ can be explicitly defined.
\end{corollary}
A related result was obtained by Krachun and the second author in \cite{KP26} in the particular case of the hexagonal lattice, where it was shown that the displacement of a self-avoiding walk of length $n$ is $O(n/\log n)$.\footnote{This bound on the displacement was allegedly improved by OpenAI to $n^{3/4+o(1)}$ in the case of the hexagonal lattice, as part of the repository \url{https://github.com/openai/math} (cf.\ group 237) released in October 2026.} These proofs on the hexagonal lattice make strong use of the parafermionic observable, and they therefore do not apply to other lattices. Our result is the first quantitative improvement of the result of \cite{DCHammond2013} on $\mathbb{Z}^d$. 

As an application of Theorem 1.1 in \cite{DCHammond2013}, Hutchcroft obtained a non-quantitative improvement to the Hammerlsey-Welsh bound $c_n\leq e^{o(\sqrt{n})}\mu_c^n$ as $n\to \infty$ \cite{Hutchcroft2018}. Using our Theorem~\ref{thm: main theorem}, we obtain the improved bound $c_n\leq e^{\Psi(n)}\mu_c^n$ for some explicit function $\Psi$ such that $\Psi(n)/\sqrt{n}\to 0$ as $n\to\infty$. Computing the precise asymptotic of $\Psi$ is beyond the scope of the current paper. 

The key ingredient in the main result of \cite{DCHammond2013} is the proof that the Ces\'aro limit of $b_k \mu_c^{-k}$ is equal to $0$, where $b_n$ denotes the number of bridges of length $n$ starting from $0$:
\begin{equation}\label{eq: cesaro}
\lim_{T\to\infty}\sum_{k=1}^T \frac{b_k \mu_c^{-k}}{T}=0.
\end{equation}
Along the way, we obtain the following quantitative improvement of \eqref{eq: cesaro} which we believe is of independent interest.

\begin{theorem}\label{thm: b_n}
There exists $C>0$ such that for every $T\geq 1$,
\[
\sum_{k=1}^T \frac{b_k \mu_c^{-k}}{T}\leq \frac{C}{\log(\log(T))}.
\]
\end{theorem}

Theorem~\ref{thm: b_n} tells us that on average we have $b_n\leq C \mu_c^n/\log(\log(n))$, thus improving on the standard bound $b_n\leq \mu_c^n$. Similar bridge partition functions have also been studied in the case of the hexagonal lattice, where it was shown in \cite{BBDCG14} that the critical partition function of bridges of height $T$ tends to $0$ (at a non-quantitative rate), and later improved to logarithmic decay in \cite{GM20} and polynomial decay in \cite{KP26}. The latter fact was used in \cite{KP26} together with the summability of the critical boundary $2$-point function in the upper half-plane to show that on the hexagonal lattice the typical displacement is at most $O(n/\log n)$. The latter is only known for the hexagonal lattice among low-dimensional lattices and its proof uses the integrability features of the model on this particular graph. Let us also remark that the optimal decay of $b_n \mu_c^{-n}$ is expected to be polynomial, and under this assumption, our proof gives the bound $O(n/(\log n)^{\varepsilon})$ on the displacement for some $\varepsilon>0$; see Remark \ref{rem:optimal_T}.

\subsection{Overview of the proof} 

On a high-level, this paper can be viewed as a proof-theoretic analysis of \cite{DCHammond2013}, in the sense that we break down the proof of Theorem \ref{res:main} given in \cite{DCHammond2013} into several different components, provide a specific quantitative version of each of these whose relationship to its qualitative counterpart can be formally explained using ideas from logic, and then reassemble them to arrive at Theorem \ref{thm: main theorem}, where the precise quantitative form of each component is crucial in our ability to reassemble the pieces. Part of the significance of the proof-theoretic aspect of the paper is that it represents a methodology that can be applied systematically to obtain other quantitative results in probability theory, building on a body of work on proof theory in probability that has been established over the last few years. To the authors' knowledge, this is the only time methods from logic have been applied to self-avoiding random walks. For that reason, this connection to proof theory is described in detail in Section \ref{sec:structure} below. However, while logic provides the overall structure, our proof can be formulated and understood without reference to it, and indeed, a number of crucial \emph{non-logical} innovations were required, which we now discuss. 

In the first part of the proof (Section \ref{sec:a}), we obtain an estimate on the decay of 
\[
\psaw{n}\left(\max\left\{\norm{\gamma_k}\mid 0\leq k\leq n\right\}\geq vn\right)
\]
in terms of an arbitrary function $h$ which encodes the decay of $\sum_{k=1}^T \frac{b_k \mu_c^{-k}}{T}$ to $0$. To this end, we gather two key observations. First, it suffices to control the height of a typical bridge of length $n$, due to the Hammersley-Welsh bridge decomposition operation. Second, as observed by Kesten, bridges have a natural renewal structure, and Theorem~\ref{thm: b_n} and standard large-deviations estimates imply an explicit fast decaying bound on the probability a bridge of length $n$ has many renewal points. To show that it is unlikely for a bridge $\gamma$ to have large height, we apply successive unfolding operations based on the zigzags of \cite{DCHammond2013} until we end up with a bridge with a large number of renewal points. The main idea for this part of our proof is using the fast decay of the probability of having many renewal points to control the cost of unfolding at each step of the iteration. Then the bound $O(n/\exp(c\log^* n))$ on the displacement arises from the iterative nature of this approach applied to a concrete function $h$, and the need to choose the height of the original bridge $\gamma$ to be large enough for the iteration to close.

The second part consists of establishing Theorem~\ref{thm: b_n} and the provision of the concrete rate of decay $h$, whose proof we now describe. The renewal structure of bridges allows us to study them by decomposing them into irreducible components called \emph{irreducible bridges} sampled under the probability measure $\pisab$ defined by $\pisab(\gamma)=\mu_c^{-|\gamma|}$, where $|\gamma|$ denotes its length. This measure can be used to define a probability measure $\pisabn$ on bridges of infinite length by concatenating i.i.d.\ irreducible bridges. In \cite{DCHammond2013}, it was shown by contradiction that 
\begin{equation}\label{eq: infinite}
\eisab[|\gamma|]=\infty.
\end{equation}
In this paper, we prove the following quantitative improvement, which is a reformulation of Theorem~\ref{thm: b_n} and whose proof is presented in Section~\ref{sec:b}. In what follows, we write $|\gamma|^{\leq M}:=|\gamma|\mathbf{1}_{|\gamma|\leq M}$.

\begin{theorem}
\label{thm:stepb}
For all $v\geq 1$ we have $\eisab\left[|\gamma|^{\leq h(v)}\right]>v$ where $h(v)=e^{e^{O(v)}}$.
\end{theorem}

The original argument of \cite{DCHammond2013} for proving \eqref{eq: infinite} goes as follows. On the one hand, assuming for a contradiction that $\eisab[|\gamma|]<\infty$, one can show using the law of large numbers that a typical infinite bridge sampled from $\pisabn$ has positive speed which is almost surely constant. On the other hand, by the law of large numbers and an ergodicity argument, a typical bridge sampled from $\pisabn$ has a positive density of a certain class of renewal points called \emph{diamond points}. Applying a geometric construction involving the so-called stickbreak operation and the multi-valued map principle, one obtains bridges that have different speed than the original one. Hence it follows that under $\pisabn$, the speed cannot be constant with probability $1$, thus leading to a contradiction. 

The argument described above is fundamentally non-constructive. However, the proof theoretic approach (as detailed in Section \ref{sec:structure}) allows us to deal with nonconstructive reasoning in a computational way, and in addition, a key innovation of this paper is a proof strategy that completely avoids the use of the law of large numbers and of ergodicity arguments. Interestingly, the presence of the law of large numbers was not in itself an obstacle: an earlier version of our proof involved a more sophisticated proof-theoretic construction which replaced instances of the law of large numbers with a finitary counterpart similar in style to that for pairwise independent random variables established in \cite{Neri2026}. This however resulted in a version of Theorem \ref{thm:stepb} where $h$ had tower exponential size, as opposed to the double exponential bound obtained via the more streamlined argument provided in Section \ref{sec:b}. 

The proof of Theorem~\ref{thm:stepb} is split into two steps both involving the event $\{0\in C_n\}$ that the concatenation of $n$ irreducible bridges stays in a cone of angle $\pi/4$. First, under the assumption $\eisab\left[|\gamma|^{\leq h(v)}\right]\leq v$ for some $v$ large enough, we obtain a polynomial \emph{lower} bound for the probability of the event $\{0\in C_n\}$. Compared to the argument of \cite{DCHammond2013}, this step avoids the use of the law of large numbers and uses instead a reflection argument to keep the bridge inside the cone each time it exits it. The assumed upper bound  $\eisab\left[|\gamma|^{\leq h(v)}\right]\leq v$ is used to control the cost of each reflection and the number of times the bridge tends to exit the cone. Second, under the assumption $\eisab\left[|\gamma|^{\leq h(v)}\right]\leq v$, we obtain a similar polynomial \emph{upper} bound for the probability of the event $\{0\in C_n\}$. To obtain this bound, we follow the argument of \cite{DCHammond2013} using a stick-break operation involving a multi-valued map argument and our novelty comes in obtaining a tighter control in the number of pre-images via a geometric construction. Since both bounds on the probability of $\{0\in C_n\}$ cannot hold at the same time, the assumption that led to them cannot hold, which gives Theorem~\ref{thm:stepb}. To make this approach work, we prove a quantitative upper bound on the tail of an irreducible bridge using the Hammersley-Welsh bound---see Lemma~\ref{lem:markov}. This result is one of the two sources of exponential growth in the definition of the function $h$, with the other source arising from the repeated application of reflections.

\section{Preliminaries}

In this section, we formally define our model and gather some definitions and results that will be useful throughout the paper.

\subsection{Notation and definition of the model}

We work on the hypercubic lattice $\mathbb{Z}^d$, $d\geq 2$, where two vertices $u,v$ are connected by an edge if $\lVert u-v \rVert_1=1$. Given a vertex $v\in \mathbb{Z}^d$, we write $x(v)$ and $y(v)$ to denote the first and second coordinates of $v$, respectively. 

We write $\N$ to denote the set of non-negative integers $\{0,1,2,3,\dots\}$ and $\N^* = \N\backslash \{0\}$. We set $[a; b] = \{a,a+1,\dots,b\}$. The cardinality of a set $A$ is denoted by $|A|$.

For $n\in \N$, a \textit{walk} of length $n$ is a sequence $(\gamma_0,
\gamma_1,\ldots, \gamma_n)$ of vertices in $\mathbb{Z}^d$ such that for every $i\in [0;n-1]$, $\gamma_i$ and $\gamma_{i+1}$ are neighbours. A \textit{self-avoiding walk} is a walk, all vertices of which are 
distinct. We denote by $\SAW_n$ the set of self-avoiding walks of length $n$ starting at $0$ and by $c_n$ the cardinality of $\SAW_n$.

The \emph{connective constant} $\mu_c=\mu_c(\mathbb{Z}^d)$ is defined as
\[
\mu_c:= \lim_{n\to \infty} c_n^{1/n},
\]
where the limit exists by a standard sub-additivity argument using Fekete's lemma, see e.g.\ \cite{Poor}. The same sub-additivity argument gives that 
\begin{equation}\label{eq: c_n}
c_n\geq \mu_c^n, \qquad \forall n\geq 0.
\end{equation}
For $z\geq0$, we define the \emph{susceptibility} as $\chi(z)=\sum_{n\geq 0} c_n z^n$, and we note that $\chi(z)<\infty$ for every $z\in [0,1/\mu_c)$.

For $m,n \in \N$, let $\gamma$ and $\Tilde{\gamma}$ be two walks of lengths $m$ and $n$, neither of which need to start at $0$. The concatenation $\gamma \circ \Tilde{\gamma}$ of $\gamma$ and $\Tilde{\gamma}$ is given by
$$
(\gamma \circ \Tilde{\gamma})_k = \begin{cases}
    \gamma_k &k \leq m\\
    \gamma_m + (\Tilde{\gamma}_{k-m}-\Tilde{\gamma}_0) &m+1 \leq k \leq m + n.
\end{cases}
$$

\subsection{Bridges and renewal theory}

A \textit{self-avoiding bridge} $\gamma$ of length $n \in \N$ on $\mathbb{Z}^d$ is a self-avoiding walk of length $n$ such that $y(\gamma_0) < y(\gamma_i) \leq y(\gamma_n)$ for all $i \in [1;n]$. The \textit{span} of $\gamma$ is defined as
$$
\spn(\gamma) := y(\gamma_n) - y(\gamma_0).
$$
Let $\SAB_n$ denote the set of self-avoiding bridges of length $n$ that start at $0$, $b_n$ the cardinality of $\SAB_n$ and $B(z)=\sum_{n\geq 0} b_n z^n$.
Let us also write
$\P_{\SAB_n}$ for the uniform law on $\SAB_n$.

Let $\SAW = \bigcup_{n \geq 0} \SAW_n$ and $\SAB = \bigcup_{n \geq 0} \SAB_n$. When clear from context, we will occasionally omit the term ``self-avoiding'' when referring to walks and bridges.

We now recall some important properties and results regarding bridges. The first statement we recall is due to Hammersley and Welsh \cite{HammersleyWelsh} and relates the number of bridges to the number of self-avoiding walks.

\begin{proposition}\label{HW}
Let $d\geq 2$. There exists a constant $C_{HW} > 0$ such that, for each $n \in \N$,
\begin{equation}\label{HW_2}
e^{-C_{HW}\sqrt{n}}c_n \leq b_n \leq c_n.
\end{equation}
\end{proposition}
\noindent We will refer to the Hammersley-Welsh constant as $C_{HW}$ throughout the paper. 

The Hammersley-Welsh bridge decomposition that gives the above result allows to deduce the following result that relates $\chi(z)$ and $B(z)$. See \cite[Corollary 2.8]{LecturesSAW} for a proof.

\begin{proposition}\label{prop:chi-B bound}
For every $z\in [0, 1/\mu_c)$ we have
$$
\chi(z)\leq \frac{1}{z}e^{2(B(z)-1)}.
$$
\end{proposition}

An important property of bridges is that any bridge can be decomposed into irreducible pieces that we call \textit{irreducible bridges}. This observation allows for the development of a renewal theory for bridges (see \cite{KestenI}, \cite{KestenII}). The set $R_\gamma$ of \textit{renewal points} of a $\gamma \in \SAB_n$ is the set of points $\gamma_i$, with $i \in [0; n-1]$, for which both $(\gamma_0,\ldots,\gamma_i)$ and $(\gamma_i,\ldots, \gamma_n)$ are bridges. We write $r_k$ to denote the index of the $k$th renewal point. An \emph{irreducible} bridge is a bridge $\gamma$ of length $n\geq 1$ such that $\gamma_k$ is not a renewal point for any $k \in [1;n-1]$. Here we emphasise that irreducible bridges are not allowed to have zero length. 

\begin{figure}[ht]
\centering
\begin{tikzpicture}[x=0.34cm,y=0.32cm]
\draw[walk] 
plot coordinates {(0,0) (0,1) (2,1) (2,3) (1,3) (1,2) (-2,2) (-2,1) (-4,1) (-4,2) (-3,2) (-3,3) (-1,3) (-1,4) (4,4) (4,2) (5,2) (5,5) (4,5) (4,6) (2,6) (2,5) (0,5) (0,7) (-1,7) (-1,8) (-2,8) (-2,9) (1,9) (1,10) (2,10) (2,11) (-4,11) (-4,10) (-7,10) (-7,12) (-5,12) (-5,16) (-3,16) (-3,15) (-2,15) (-2,13) (-4,13) (-4,12) (1,12) (1,15) (-1,15) (-1,16)};

\foreach \x/\y in {0/0,5/4,0/6,-1/7,-2/8,1/9,-7/11, -1/16}{\draw[cone] (-7.5,\y)--(5.5,\y);\node[point] at (\x,\y) {};}
\end{tikzpicture}
\caption{An example of a bridge and its decomposition into irreducible bridges. The renewal points are highlighted in blue.}
\label{fig: renewal points}
\end{figure}

Let $\iSAB$ be the set of irreducible bridges of arbitrary length whose initial vertex is $0$, and for $z\geq 0$, define $B_{\textup{irr}}(z)=\sum_{\gamma \in \iSAB} z^{|\gamma|}$. Every bridge can be uniquely decomposed into a finite number of irreducible bridges. This implies that for every $z\in [0,1/\mu_c)$,
\begin{equation}\label{eq:decomposition}
B(z)=\sum_{k\geq 0} (B_{\textup{irr}}(z))^k=\frac{1}{1-B_{\textup{irr}}(z)}.
\end{equation}
Since $B(1/\mu_c)=\infty$, which follows from Proposition~\ref{prop:chi-B bound} and \eqref{eq: c_n}, we obtain the following result, originally due to Kesten \cite{KestenI}.

\begin{lemma}[Kesten]\label{kesten}
We have $B_{\textup{irr}}(1/\mu_c)=\sum_{\gamma \in \iSAB} \mu_c^{-|\gamma|} = 1$.
\end{lemma}

We define the measure $\P_\iSAB$ on $\iSAB$ by setting \[
\P_\iSAB(\gamma) = \mu_c^{-|\gamma|}, \qquad \gamma\in \iSAB.
\]
By Lemma \ref{kesten}, $\P_\iSAB$ is a probability measure. The following tail bound on the lengths of irreducible bridges will be crucial in what follows.

\begin{lemma}
\label{lem:markov}
Let $d\geq 2$. There exists $C=C(d)>0$ such that for every $k\geq 2$ we have
\[
\pisab\left(|\gamma|\geq k\right)\leq \frac{C}{\log(k)}.
\]
\end{lemma}
\begin{proof}
It suffices to prove the result for every $k>1+\mu_c/e^2$. To this end, recall \eqref{eq: c_n} and the definition of $\chi(z)$. Then, for every $z\in [0,1/\mu_c)$, $\chi(z)\geq \frac{1}{1-z\mu_c}$. For $z_k=\frac{1-1/k}{\mu_c}$ this gives $\chi(z_k)\geq k$. Using Proposition~\ref{prop:chi-B bound} we obtain
$$
B(z_k)\geq \frac{\log(z_k \chi(z_k))}{2}+1\geq \frac{\log(k-1)-\log(\mu_c)+2}{2}.
$$
Combined with \eqref{eq:decomposition} this, in turn, implies that 
\begin{equation}\label{eq:lower bound}
B_{\textup{irr}}(z_k)=1-\frac{1}{B(z_k)}\geq 1-\frac{2}{\log(k-1)-\log(\mu_c)+2},    
\end{equation}
where here we note that because of our choice of $k$, we have $\log(k-1)-\log(\mu_c)+2>0$.

On the other hand, using that $(1-1/k)^n\leq 1$ for all $n\geq 1$, and $(1-1/k)^n\leq 1/e$ for every $n\geq k$, we get 
\begin{equation}\label{eq:upper bound}
B_{\textup{irr}}(z_k)\leq \mathbb{P}_{\textup{iSAB}}(|\gamma|< k)+\frac{1}{e} \mathbb{P}_{\textup{iSAB}}(|\gamma|\geq k)=1-(1-1/e)\mathbb{P}_{\textup{iSAB}}(|\gamma|\geq k).    
\end{equation}
Combining \eqref{eq:lower bound} with \eqref{eq:upper bound}, we obtain
$$
\mathbb{P}_{\textup{iSAB}}(|\gamma|\geq k)\leq \frac{2e}{(e-1)(\log(k-1)-\log(\mu_c)+2)}.
$$
The desired bound follows from choosing $C>0$ to be large enough.
\end{proof}

It will be useful in our analysis to define the function $\phi:(0,1]\rightarrow \mathbb{N}$ by 
\begin{equation}\label{eq: phi def}
\phi(\lambda)=\min\big\{k\geq 1/\lambda: \, \mathbb{P}_{\textup{iSAB}}(|\gamma|\geq k)\leq \lambda\big\}.
\end{equation}
Note that by definition, $\phi(\lambda)\geq 1/\lambda$ and by Lemma~\ref{lem:markov}, $\phi(\lambda)\leq e^{C/\lambda}$.

We construct a probability measure on infinite self-avoiding bridges as the product measure $\infP$. That means, we consider $\gamma: \N \to \Z^d$ an infinite self-avoiding walk distributed as the concatenation of an infinite number of irreducible self-avoiding bridges with distribution $\irrP$. We extend the notion of renewal point to this infinite bridge in the natural way.

We recall the following lemma on the law of the infinite half-plane walk conditioned on the existence of renewal points.

\begin{lemma}[Lemma 3.8 of \cite{DCHammond2013}]\label{isab_renewal} For each $n \in \N$, we have that
    $$
    b_n \mu_c^{-n} = \infP[\gamma_n \in \mathsf{R}_\gamma].
    $$
    The conditional distribution of $\gamma_{[0;n]}$ under $\infP[\ \cdot \mid \gamma_n \in \mathsf{R}_\gamma]$ coincides with $\P_{\SAB_n}$.
\end{lemma}

Finally, we note that a key tool in the study of self-avoiding walk, and in this article in particular, is the multi-valued map principle, which we now state.
Let $f$ be a multi-valued map from $A$ to $B$. Write $f^{-1}(b) := \{a \in A : b \in f(a)\}$. Then
$$
|A| \leq \frac{\max_{b \in B}|f^{-1}(b)|}{\min_{a \in A}|f(a)|} \cdot |B|.
$$

\section{Connections to proof theory}
\label{sec:structure}

Our approach was partly inspired by the \emph{proof mining} programme, which, using ideas and techniques from mathematical logic (and proof theory in particular) seeks to answer the following question: by analysing an existing proof from a logical perspective, can we obtain new insights and information? Over the last 35 years or so, a large body of work has formed in which new information has been systematically extracted from proofs, which very often takes the form of \emph{explicit bounds} for existential statements. Success here often involves combining proof-theoretic techniques with ordinary mathematical reasoning, and our route to the main result of the present paper involved studying the original, qualitative proof of \cite[Theorem 1.1]{DCHammond2013} using G\"odel's Dialectica interpretation \cite{Goedel1958} in a form in which it has recently been applied in probability theory (see \cite{NeriPischke2024} for a logical account and \cite{NeriPowell2025a} for a recent case study).

A detailed account of proof mining can be found in the monograph \cite{Kohlenbach2008} as well as the surveys \cite{Kohlenbach2019,KohlenbachOliva2003}. The purpose of this section is simply to highlight, very informally, how the construction of our main bound can be viewed from a proof-theoretic perspective. Not only do we consider this connection between logic and probability to be of interest in its own right, but it also provides a convenient map of the technical results that follows. Discussions around logic are restricted to this section, and the remainder of the paper can be understood without reference to any of the theory we allude to here (a fact which also justifies our informal presentation). Indeed, one of the virtues of proof mining is that once a bound has been obtained, any logical apparatus can normally be suppressed and the quantitative result presented in ordinary mathematical language.

\subsection{The Dialectica interpretation} We start with a brief sketch of G\"odel's Dialectica interpretation (see \cite{AvigadFeferman1998} or \cite{Kohlenbach2008} for a modern account). Roughly speaking, this takes a statement $A$ in some formal system, and interprets it as an equivalent statement of the form $A^D=\exists x\forall y A_D(x,y)$, where $x,y$ are (potentially empty) sequences of variables, and $A_D(x,y)$ is quantifier-free and thus, in some sense, ``computationally neutral''. This translation on formulas is defined inductively over their logical structure, the most interesting and relevant clause being implication. Here, if $A^D=\exists x\forall y A_D(x,y)$ and $B^D=\exists u\forall v B_D(u,v)$, the Dialectica characterises the computational content of the implication $A\Rightarrow B$ as a pair of functions $f(x)$ and $g(x,v)$ such that
\begin{equation}
\forall x,v\left(A_D(x,g(x,v))\implies B_D(f(x),v)\right).\label{dial}\tag{$\ast$}
\end{equation}
We now arrive at the important question: given a proof of $A$, can we produce (a bound on) some explicit term $t$ satisfying $\forall y A_D(t,y)$. In the case that we have proved an implication $A\Rightarrow B$, this becomes the question of whether we can produce explicit (bounds on) functions\footnote{The precise way in which the notion of bounds extends to functions and functionals is formally dealt with via the concept of majorization and the monotone variant of the Dialectica interpretation \cite{Kohlenbach1996}, though we do not mention any details here.} $f(x)$ and $g(x,v)$ satisfying \eqref{dial}.

In order to answer this question rigorously, we must first fix a formal system for representing the kind of proofs we want to analyse, and this alone is an extremely delicate matter, particularly in probability theory. Here, a naive formalisation in which basic concepts such as probability measures and random variables are formally encoded in terms of functions on natural numbers leads to an unwieldly system of vast complexity. However, the last few years have seen the development of specialised systems for probability \cite{NeriOlivaPischke2026,NeriPischke2024} (and a forthcoming paper by the last two authors with Neri on the treatment of unbounded random variables) along with so-called \emph{metatheorems} which guarantee that explicit bounds can be extracted from large classes of proofs in probability theory. Moreover, these bounds satisfy striking uniformities. A key insight, made formal in the aforementioned works, is that one can make computational sense of probabilistic statements involving infinite unions or intersection by replacing those events with \emph{finitary approximations} that are formulated using only finitely many Boolean operations on atomic events, where the latter include events of the form $X\geq a$ for $X$ a random variable.

Under these assumptions, the three key steps in the present paper (Theorems \ref{thm:stepa_inverse}, \ref{thm:contradiction-a} and \ref{thm:contradiction-b}) can be viewed as the construction of (bounds on) functions for the Dialectica interpretation of a series of implicative statements (\eqref{imp1}, \eqref{imp2} and \eqref{imp3} below), which together reflect the broad structure of the original proof of Theorem \ref{res:main} in \cite{DCHammond2013}. Crucially, these implications are individually proven \emph{constructively} (that is, without using proof by contradiction), and the functions arise naturally from those proofs. However, as we will see, combining the interpretations of these implications allows us to give a computational meaning to a global strategy that \emph{does} use contradiction.

\subsection{The logical structure of the proof in \cite{DCHammond2013}} We first turn to the logical structure of Theorem \ref{res:main}, which is
\[
T=\forall v\in \QQ_+\, \exists \varepsilon\in \QQ_+ \exists n_0\in\NN\, \forall n\geq n_0\, T_D(v,\varepsilon,n)
\]
for
\[
T_D(v,\varepsilon,n):=\psaw{n}\left(\max\left\{\norm{\gamma_k}\mid 0\leq k\leq n\right\}\geq v n\right)\leq e^{-\varepsilon n},
\]
noting that the event $A$ inside the probability is a finite union of events of the form $\norm{\gamma_k}\geq vn$, in line with our earlier observation.\footnote{We note that, technically, there is an additional quantifier implicit in the inequality $\PP(A)\leq e^{-\varepsilon n}$ between real numbers, but for logical reasons that we do not expand on here, this can be viewed as benign.} Here we prefer to represent the positive reals $(0,\infty)$ in a finite way via the positive rationals $\QQ_+$, and note that all functions $(0,\infty)\to (0,\infty)$ produced in this paper are in particular computable in the sense of computability theory when restricted to positive rationals. 

The proof in \cite{DCHammond2013} can be divided into two parts. First, one shows that if the conclusion of Theorem \ref{res:main} is violated, then $\eisab[|\gamma|]<\infty$ (\cite[Theorem 2.3, Corollary 2.4]{DCHammond2013}), or equivalently:
\begin{equation}
\eisab[|\gamma|]=\infty\implies T.\label{imp1}\tag{I1}
\end{equation}
One then shows that $\EE_{\isab}(|\gamma|)=\infty$ (\cite[Theorem 2.5]{DCHammond2013}) by contradiction, and we give a modified proof of the latter that can be decomposed along the lines of \cite{DCHammond2013}. Specifically, letting $\{0\in C\}$ denote the event that the concatenation of infinitely many independent samples of $\pisab$ stays in a cone of angle $\pi/4$, one can simultaneously show that
\begin{equation}
\eisab[|\gamma|]<\infty\implies \pisabn(0\in C)>0\label{imp2}\tag{I2}
\end{equation}
(which forms a simplified analogue of \cite[Proposition 4.2]{DCHammond2013}), but also
\begin{equation}
\eisab[|\gamma|]<\infty\implies \pisabn(0\in C)=0,\label{imp3}\tag{I3}
\end{equation}
the latter via a modified version of the stickbreak operation given in \cite[pp.~22--25]{DCHammond2013}. An important improvement in our case is that both of these individual implications can be proven constructively, with the law of large numbers and ergodicity argument in \cite{DCHammond2013} eliminated in favour of a strategy based on Levy's maximal inequality and a Borel-Cantelli type argument. Our overall proof that $\EE_{\isab}(|\gamma|)=\infty$ is still by contradiction, but here the Dialectica interpretation of the corresponding implications allows us to deal with this computationally, as we now show.

\subsection{The Dialectica interpretation of \eqref{imp1}, \eqref{imp2} and \eqref{imp3}}

We now sketch the Dialectica interpretation of the statements we have identified. These not only provide insight into the structure of the key computational theorems of the paper, but also how they combine to prove the main result (Theorem \ref{thm: main theorem}). We first consider \eqref{imp1}, for which we need to understand the quantifier structure of $\eisab[|\gamma|]=\infty$. In particular, since we \emph{cannot} assume that $|\gamma|\in L^1$, we instead consider truncations $|\gamma|^{\leq M}:=|\gamma|\mathbf{1}_{|\gamma|\leq M}$ and formulate $\eisab[|\gamma|]=\infty$ as $\forall v\in \QQ_+\, \exists M\in \QQ_+\, \eisab[|\gamma|^{\leq M}]>v$. Then it is easy to see that \eqref{imp1} is equivalent to
\[
\exists h\, \forall v\in \QQ_+\left(\eisab[|\gamma|^{\leq h(v)}]>v\right)\implies \exists \varepsilon,n_0\, \forall \delta\in \QQ_+\forall n\geq n_0(\delta)\, T_0(\delta,\varepsilon(\delta),n).
\]
where now $h,\varepsilon:\QQ_+\to \QQ_+$ and $n_0:\QQ_+\to \NN$ are functions. Because we are later able to provide an explicit function $h$ satisfying the premise of the above, we actually only require a pair of functionals $F_1(h,\delta),F_2(h,\delta)$ which satisfy the `forward' direction of the Dialectica interpretation i.e.\ witnessing $\varepsilon$ and $n_0$, and we do not need to identify the `backward' functional witnessing $\forall v$ in the premise. In other words, it is sufficient to find functionals $F_1,F_2$ such that
\[
\forall h\, \left(\forall v\in \QQ_+\left(\eisab[|\gamma|^{\leq h(v)}]>v\right)\implies \forall \delta\in \QQ_+,\forall n\geq F_2(h,\delta)\, T_0(\delta,F_1(h,\delta),n)\right).
\]
The explicit construction of such an $F_1,F_2$ is directly implied by Theorem \ref{thm:stepa_inverse}, as best seen in \eqref{eq:bound_delta_alpha}. In particular, these functionals allow us to transfer an asymptotic bound for any given $h$ to one for the functions $F_1(h,\cdot)$ and $F_2(h,\cdot)$. In this paper, we are able to show that $h(v)\in e^{e^{O(v)}}$, and our main bound in Theorem \ref{thm: main theorem} then follows from a complexity analysis of the functionals $F_1,F_2$ in this specific case.

Our construction of a concrete $h(v)\in e^{e^{O(v)}}$ is, in turn, obtained by combining the Dialectica interpretations of \eqref{imp2} and \eqref{imp3}. For this, we need to understand the quantifier structure of $\pisabn(0\in C)$. Here, we observe that 
\[
\{0\in C\}=\bigcap_{n\in\NN}\{0\in C_n\}
\]
where $\{0\in C_n\}$ denotes the event that the concatenation of $n$ independent irreducible bridges stays in the cone. Since the sequence $\{0\in C_n\}$ is monotone decreasing, the full quantifier structure of the statement $\pisabn(0\in C)>0$ becomes $\exists \delta\in \QQ_+\, \forall n\in\NN\left(\pisabn(0\in C_n)\geq\delta\right)$, where now we observe that $\{0\in C_n\}$ can be written as a finite intersection of basic events of the form $(y+x)(\gamma_k)\geq 0$ or $(y-x)(\gamma_k)\geq 0$. Thus \eqref{imp2} is equivalent to
\[
\exists v\in \QQ_+\, \forall M\in \QQ_+\left(\eisab[|\gamma|^{\leq M}]\leq v\right)\implies \exists \delta\in \QQ_+\, \forall n\in\NN\left(\pisabn(0\in C_n)\geq\delta\right)
\]
and its Dialectica interpretation is given by a pair of functions $f(v)$ and $g(v,n)$ satisfying
\begin{equation}
\forall v\in \QQ_+,n\in\NN\left(\eisab[|\gamma|^{\leq g(v,n)}]\leq v \implies \pisabn(0\in C_n)\geq f(v)\right).\label{dial1}\tag{D1}
\end{equation}
The construction of such functions is given in Theorem \ref{thm:contradiction-a}. In an analogous way, we can show that \eqref{imp3} is equivalent to
\[
\exists v\in \QQ_+\, \forall M\in \QQ_+\left(\eisab[|\gamma|^{\leq M}]\leq v\right)\implies \exists n\, \forall \delta\in \QQ_+\left(\pisabn(0\in C_{n(\delta)})<\delta\right)
\]
for $n:\QQ_+\to \NN$ a function, whose Dialectica interpretation is now a pair of functions $N(v,\cdot)$ and $L(v,\delta)$ satisfying
\begin{equation}
\forall v,\delta\in \QQ_+\left(\eisab[|\gamma|^{\leq L(v,\delta)}]\leq v\implies \pisabn(0\in C_{N(v,\delta)})<\delta\right).\label{dial2}\tag{D2}
\end{equation}
The construction of such functions is given in Theorem \ref{thm:contradiction-b}. Fixing $v>0$, we now combine \eqref{dial1} and \eqref{dial2} by setting $\delta:=f(v)$ in \eqref{dial2}, $n:=N(v,\delta)$ in \eqref{dial1}, and defining
\[
h(v):=\max\left\{g(v,N(v,f(v))),L(v,f(v))\right\},
\]
we obtain
\[
\forall v\left(\eisab\left[|\gamma|^{\leq h(v)}\right]\leq v \text{ is false}\right)
\]
and therefore $\forall v\left(\eisab\left[|\gamma|^{\leq h(v)}\right]> v\right)$. A complexity analysis of the functions that make up $h$ as defined in Theorems \ref{thm:contradiction-a} and \ref{thm:contradiction-b} gives us our asymptotic bound $h(v)=e^{e^{O(v)}}$.

\subsection{Summarising the role of proof theory}
Logical insights were used in the first instance to help decompose the proof in \cite{DCHammond2013}, and identify, via the Dialectica interpretation, the correct quantitative formulations of the key Lemmas \eqref{imp1}, \eqref{imp2} and \eqref{imp3}. Numerous mathematical (as opposed to logical) insights and simplifications were then needed to find witnesses for these quantitative reformulations. However, our ability to recombine these individual quantitative results to arrive at the overall bound, and in particular eliminate the use of proof by contradiction in the original proof, was only possible due to the precise way in which the main lemmas were reformulated. Because proof theory was only used in a specific, local way, our entire argument can be presented in ordinary mathematical language, and without reference to any of the underlying proof-theoretic techniques.

\section{Explicit bound on the displacement of self-avoiding walk}
\label{sec:a}

This section is dedicated to the proof of Theorem \ref{thm: main theorem} under the assumption that Theorem \ref{thm:stepb} holds. We will mainly work with bridges. First, we show in Section \ref{sec:a_asymptotics} that, under the assumption that irreducible bridges are typically long under $\irrP$, one has that finite bridges of given length cannot have many renewal points. Then, in Section \ref{sec:a_unfolding}, we use an iteration procedure similar to the one of \cite[Proposition 3.2]{DCHammond2013} to show that tall bridges are uncommon under the same assumption. Finally, in Section \ref{sec:a_analysis} we use that bound to deduce Theorem \ref{thm: main theorem}.

\subsection{A concentration result for the number of renewal points} \label{sec:a_asymptotics}

Define 
\begin{equation*}
q(M) := \irrE[|\gamma| \land M].    
\end{equation*}
We now introduce a useful concentration result on the number of renewal points, recalling that $C_{HW}$ is the Hammersley-Welsh constant:

\begin{lemma}\label{lemma:bound_renewal}
For every $n \in \N^*$ and $M >0$,
$$
\P_{\SAB_n}\left(|R_\gamma| \geq \frac{2n}{q(M)} \right) \leq \exp\left(-\frac{n}{5M} + C_{HW} \sqrt{n}\right).
$$
\end{lemma}

\begin{proof}
Throughout this proof, let $h(x)= (1+x)\log(1+x) -x$. Consider a sequence of i.i.d.\ irreducible bridges $\gamma^i$ sampled from $\pisab$. For all $i\in \N^*$, define $X_i = \infE[|\gamma^{i}| \land M] - |\gamma^{i}| \land M$. Notice that the random variables $X_i$ are independent, $\infE[X_i]=0$, and $|X_i| \leq M$ for all $i$. Then, by Bennett's inequality \cite[Theorem 2.9]{BLM2013}, for every $k\geq 1$ and every $t > 0$,
$$
\infP\left(\sum_{i=1}^k X_i \geq t \right) \leq \exp\left(-\frac{\sigma^2}{M^2} h \left(\frac{M t}{\sigma^2}\right)\right),
$$
where $\sigma^2 = \sum_{i=1}^k \infE[X_i^2]$.
Notice that the expression $\frac{\sigma^2}{M^2} h \left(\frac{M t}{\sigma^2}\right)$ is decreasing in $\sigma^2$, hence, using that $\infE[X_i^2] \leq M \, \infE[ |\gamma^{i}| \land M] = M q(M)$, we get that
$$
\infP\left(\sum_{i=1}^k (|\gamma^{i}| \land M )\leq -t + k q(M)  \right) \leq \exp\left(-\frac{k q(M)}{M} h \left(\frac{t}{k q(M)}\right)\right).
$$
In particular, for $t = \lambda k  q(M)$ with $0 < \lambda < 1$,
$$
\infP \left(\sum_{i=1}^k (|\gamma^{i}| \land M) \leq (1-\lambda) k q(M) \right) \leq \exp\left(-\frac{k q(M)}{M} h \left(\lambda\right)\right).
$$
Define $S_k = \sum_{i=1}^k |\gamma^{i}|$ and apply the previous inequality for $k=\left\lceil \frac{n}{(1-\lambda)q(M)} \right\rceil$ to get
$$
\infP\left(S_{\left\lceil \frac{n}{(1-\lambda)q(M)} \right\rceil} \leq n\right) \leq \exp\left(-\frac{h(\lambda)}{1-\lambda} \frac{n}{M}\right).
$$
Note that for every $n \in \N$,
$$
\infP\left(|R_\gamma[1,n]| \geq  \frac{n}{(1-\lambda)q(M)} \right) = \infP\left(S_{\left\lceil \frac{n}{(1-\lambda) q(M)} \right\rceil} \leq n\right).
$$
Proposition \ref{HW} and Lemma \ref{isab_renewal} finally give
\begin{align*}
\P_{\SAB_n}\left(|R_\gamma| \geq \frac{n}{(1-\lambda) q(M)} \right) 
&= \infP\left(|R_\gamma[0,n-1]| \geq  \frac{n}{(1-\lambda) q(M)} \, \Big| \, \gamma_n \in R_\gamma\right)\\
&= \infP\left(|R_\gamma[1,n]| \geq  \frac{n}{(1-\lambda) q(M)} \, \Big| \, \gamma_n \in R_\gamma\right)\\
&\leq \exp\left(-\frac{h(\lambda)}{1-\lambda}\frac{n}{M} + C_{HW} \sqrt{n}\right).
\end{align*}
We conclude the proof by taking $\lambda= 1/2$, for which we have $2h(1/2) > 1/5$.
\end{proof}

\subsection{Unfolding bridges}\label{sec:a_unfolding}
In this section, we will work with finite bridges. First, we will reintroduce the concept of zigzag and some elementary results as presented in \cite{DCHammond2013}. We then prove Proposition \ref{prop:unfolding_tree}, which is a quantitative analogue of \cite[Proposition 3.2]{DCHammond2013}. The bound is obtained through an iterative argument and, as we shall see, is a function of the maximal number of steps of the iteration.

The following notions of zigzags and unfolding maps are taken from \cite{DCHammond2013}. Let $\gamma \in \SAB_n$. A \textit{zigzag} of $\gamma$ is a pair $(i,j)$ of indices belonging to $[0;n]$ that satisfy $i < j$ and for which
\begin{itemize}
    \item the largest value of $k$ at which the maximum of $y(\gamma_k)$ over $1 \leq k \leq j$ is achieved equals $i$, and
    \item the largest value of $k$ at which the minimum of $y(\gamma_k)$ over $i \leq k \leq n$ is achieved equals $j$.
\end{itemize}
Let $\ZigZag_\gamma$ denote the set of zigzags of $\gamma$. For $(i,j) \in \ZigZag_\gamma$, we call the image of the subwalk $\gamma_{[i,j]}$ the \textit{central section} of $(i,j)$. The indices $i$ and $j$ are respectively called the \textit{point of zig} and the \textit{point of zag} of $(i,j)$.

Let $\mathcal{R}_v : \R^d \to \R^d$ denote the orthogonal reflection with respect to the hyperplane $\{u \in \R^d:y(u) = y(v)\}$. For $\gamma \in \SAB_n$ and $(i,j) \in \ZigZag_\gamma$, let $\operatorname{Unf}_{i,j}(\gamma) \in \SAB_n$ be given by $\operatorname{Unf}_{i,j}(\gamma) = \gamma_{[0,i]} \circ \mathcal{R}_{\gamma_i}(\gamma_{[i,j]}) \circ \gamma_{[j,n]}$.

The following lemma is straightforward and is stated without proof:

\begin{lemma}\cite[Lemmas 3.3 and 3.5]{DCHammond2013}
    Let $\gamma \in \SAB_n$.
    \begin{enumerate}
        \item The central sections of the zigzags of $\gamma$ are pairwise disjoint.
        \item Let $(i,j) \in \ZigZag_\gamma$. Then $\gamma_i$ and $\operatorname{Unf}_{i,j}(\gamma)_j$ each belong to $\mathsf{R}_{\operatorname{Unf}_{i,j}(\gamma)}$.
        \item Let $(i,j) \in \ZigZag_\gamma$. Then $\ZigZag_\gamma\backslash\{(i,j)\} \subset \ZigZag_{\operatorname{Unf}_{i,j}(\gamma)}$.
        \item Let $z_1,z_2 \in \ZigZag_\gamma$. Then $\operatorname{Unf}_{z_1}\operatorname{Unf}_{z_2} (\gamma) = \operatorname{Unf}_{z_2}\operatorname{Unf}_{z_1}(\gamma)$.
    \end{enumerate}
\end{lemma}

The last item of this lemma allows us to extend the action of $\operatorname{Unf}$ so that several zigzags are undone at once. Let $\gamma \in \SAB_n$ and $Z \subset \ZigZag_\gamma$. Set $\operatorname{Unf}_Z(\gamma) \in \SAB_n$ equal to the bridge obtained by iteratively applying the maps $\operatorname{Unf}_z$ for $z \in Z$.

We now introduce some further notation that will be useful for the rest of the section. Let $n, \ell$ be natural numbers and $k, \kappa \in [1,+\infty)$ be real numbers with $2k$ an integer. For $h \in \Z$ and $\gamma \in \SAB_n$, define the set of edges $\mathcal V_{h,h+1}$ by $\mathcal V_{h,h+1}(\gamma) := \{(\gamma_i,\gamma_{i+1}) : i \in [0;n-1], \{y(\gamma_i),y(\gamma_{i+1})\}=\{h,h+1\}\}$. We define
\begin{align*}
   &\SAB_{n,k} := \left\{ \gamma \in \SAB_n : \spn(\gamma) \geq \frac{n}{k}\right\}, \\
    &\SAB_{n,k,\kappa}^{\ell} := \left\{ \gamma \in \SAB_{n,k} : |\{h \in \Z: 1 \leq |\mathcal{V}_{h,h+1}| \leq \ell\}| \geq \frac{n}{\kappa}\right\}.
\end{align*}

Let us now gather the following remarks:
\begin{itemize}
    \item $\SAB_{n,k, \kappa}^{1} = \{\gamma \in \SAB_{n,k}:|\mathsf{R}_\gamma| \geq \frac{n}{\kappa}\}$, since for all $0 \leq i \leq n-1$, $\gamma_i \in \mathsf{R}_\gamma$ if and only if $\mathcal{V}_{y(\gamma_i),y(\gamma_i)+1} =\{(\gamma_i,\gamma_{i+1})\}$.
    \item For all $\ell \geq 1$, $\SAB_{n,k, \kappa}^{2 \ell} = \SAB_{n,k, \kappa}^{2 \ell - 1}$, because each crossed layer must be crossed an odd number of times.
    \item $\SAB_{n,k} \subset \SAB_{n,k,2k}^{2k}$, since for $\gamma \in \SAB_{n,k}$, we have 
    $$|\{h \in \Z : |\mathcal{V}_{h,h+1}| > 2 k\}| \, 2 k \leq n$$
    and 
    $$|\{h \in \Z : 1 \leq |\mathcal{V}_{h,h+1}|\}| = \spn(\gamma) \geq \frac{n}{k}.$$
\end{itemize}
Define the function $\rho: [4,\infty) \to (0,1]$ by
$$
\rho(x):=
    \frac{4}{x}\left(1 + \log \left(\frac x4\right) \right).
$$
Let $\bar q$ be defined as
\begin{equation*}
    \bar q(x) := \frac{1}{2}q\left(\frac{1}{10 \rho\left(\left(\frac{x}{5}\right)^{1/3}\right)}\right)
\end{equation*}
for $x\geq 5\cdot 8^3$ and $\bar q(x)=0$ otherwise.
For $x > 0$ and $g$ a nondecreasing function, define
\begin{equation*}
  M(x,g) := \sup\{m \in \N^*: 3^{m} \leq g^{\circ m}(x)\}  
\end{equation*}
with $M(x,g):= 0$ when the set is empty.

We now state the key proposition of this section.

\begin{proposition}\label{prop:unfolding_tree}
For any $n\in\N^*$ and $\alpha\geq 8$,
$$
\P_{\SAB_n}\left(\spn(\gamma) \geq \frac{2n}{3^{M(\bar q(5\alpha^3),\bar q)}}\right) \leq \frac{\exp\left(-\frac{\rho(\alpha)}{2} n+ C_{HW} \sqrt{n}\right)}{1-\exp\left(-\frac{\rho(\alpha)}{2} n\right)}
$$
\end{proposition}

Before getting to the proof, let us give a sketch of the argument. If a bridge has neither many renewal points nor many zigzags, then one can apply an unfolding procedure to the bridge without losing too much information on the original bridge. One can then iterate this procedure as long as not too many renewal points or zigzags are created. The main observation is that if we are able to iterate a given explicit number of times, we will necessarily have enough renewal points. Otherwise, that means that at some step of the iteration we get either many renewal points or many zigzags. In that second case, there are multiple ways to unfold a given proportion of these zigzags to create renewal points, and we do so while controlling the cost. In any case, one can always obtain a walk with many renewal points without losing too much information. One can then explicitly control the thresholds for the number of renewal points and zigzags at each step so that, no matter when we end the iteration, the bound given by Lemma \ref{lemma:bound_renewal} beats the cost of the unfolding procedure.

\begin{proof}

Let $\ell \geq 3$ be an integer and $\alpha, \beta, k, \kappa$ be positive reals. Suppose for now that $\alpha \geq 8$. Let $p \in (0,1)$ and $n \in \N$. In a similar way to \cite{DCHammond2013}, we consider the following trichotomy:

\paragraph{\textbf{Case 1: Many renewal points.}}
Suppose that
$$
\left| \left\{ \gamma \in \SAB_{n,k,\kappa}^\ell : |\mathsf{R}_\gamma| \geq \frac{n}{\beta} \right\} \right| \geq \frac{p}{2}|\SAB_{n,k,\kappa}^\ell|.
$$
Then, by the first remark above,  
\begin{equation}\label{eq:case1_result}
    |\SAB_{n,k,\beta}^1| \geq \frac{p}{2}|\SAB_{n,k,\kappa}^\ell|.
\end{equation}

\paragraph{\textbf{Case 2: Few zigzags and few renewal points.}}

Suppose that
$$
\left| \left\{ \gamma \in \SAB_{n,k,\kappa}^\ell : |\ZigZag_\gamma| \leq \frac{2n}{\alpha}, |\mathsf{R}_\gamma| \leq \frac{n}{\beta} \right\} \right| \geq (1-p)|\SAB_{n,k,\kappa}^\ell|.
$$

Let $\gamma \in  \SAB_{n,k,\kappa}^\ell$ and $\gamma' = \operatorname{Unf}_{\ZigZag_\gamma}(\gamma)$. We have that $\spn(\gamma') \geq \spn(\gamma)\geq \frac{n}{k}$.
Furthermore, if $h \in \Z$ is such that $2 \leq |\mathcal{V}_{h,h+1}(\gamma)| \leq \ell$, then the successive unfoldings performed by $\operatorname{Unf}_{\ZigZag_\gamma}$ will leave the set $\mathcal{V}_{h,h+1}$ unchanged except for changing
the value of $h$ until the first zigzag whose central section includes an edge of
$\mathcal{V}_{h,h+1}(\gamma)$ is unfolded, at which time, the set $\mathcal{V}_{h,h+1}(\gamma)$ will be split into three sets: the set of edges before the point of zig, between the point of zig and the point of zag (i.e.\@ in the central section), and after the point of zag. Each one of these sets will correspond to a set $\mathcal{V}_{h',h'+1}$ for some new height $h' \in \Z$; subsequently, the smallest of these sets has at most one third as many edges as the original set. We thus see that the number
of values of $h \in \Z$ such that $1 \leq |\mathcal{V}_{h,h+1}(\gamma')| \leq \ell/3$ is at least the number of $h \in \Z$
such that $1 \leq |\mathcal{V}_{h,h+1}(\gamma)| \leq \ell$.
Therefore, we deduce, for all $\ell \geq 3$,
$$
\operatorname{Unf}_{\ZigZag_\gamma}(\SAB_{n,k,\kappa}^\ell) \subset \SAB_{n,k,\kappa}^{\lfloor\ell/3\rfloor}.
$$

Let $\phi \in \SAB^{\lfloor\ell/3\rfloor}_{n,k,\kappa}$ and suppose that $\phi = \operatorname{Unf}_{\ZigZag_\gamma} (\gamma)$ for some unknown $\gamma$. To recover $\gamma$, one only needs to know the set of points in $\phi$ that were points of zig or zag in $\gamma$ and were unfolded by $\operatorname{Unf}_{\ZigZag_\gamma}$.
If $|\ZigZag_\gamma|\leq 2n/\alpha$, this set of points has cardinality at most $4n/\alpha$. Therefore, we have
$$
\left| \operatorname{Unf}_{\ZigZag_\gamma}^{-1}(\phi) \cap  \left\{ \gamma \in \SAB_{n,k,\kappa}^\ell : |\ZigZag_\gamma| \leq \frac{2n}{\alpha} \right\} \right| \leq \sum_{j = 0}^{\lfloor 2n/\alpha \rfloor}\binom{n}{2j}.
$$
Since we supposed $\alpha \geq 8$, we have for $0 < x < 1$,
$$
\sum_{j = 0}^{\lfloor 2n/\alpha \rfloor}\binom{n}{2j}
\leq \sum_{j = 0}^{\lfloor 4n/\alpha \rfloor}\binom{n}{j}
\leq x^{-4n/\alpha} \sum_{j = 0}^{n}\binom{n}{j}x^{j} 
= x^{-4n/\alpha} (x+1)^n \leq x^{-4n/\alpha} e^{nx},
$$
which for $x = 4/\alpha$ gives
$$
\sum_{j = 0}^{\lfloor 2n/\alpha \rfloor}\binom{n}{2j} \leq \left(\frac{e\alpha}{4}\right)^{\frac{4n}{\alpha}} = e^{\rho(\alpha)n}.
$$
Applying the multi-valued map principle, we obtain

\begin{align*}
|\SAB_{n,k,\kappa}^{\lfloor\ell/3\rfloor}| 
\geq 
e^{-\rho(\alpha) n}
\left| \left\{ \gamma \in \SAB_{n,k,\kappa}^\ell : |\ZigZag_\gamma| \leq \frac{2n}{\alpha} \right\} \right| 
\geq 
(1-p) e^{-\rho(\alpha) n} |\SAB_{n,k,\kappa}^\ell|.
\end{align*}

\paragraph{\textbf{Case 3: Many zigzags and few renewal points.}}
Suppose that
$$
|\manyZZfewR| \geq \frac{p}{2}|\SAB_{n,k,\kappa}^\ell|,
$$
where $\manyZZfewR := \{\gamma \in \SAB_{n,k,\kappa}^\ell : |\ZigZag_\gamma| \geq \frac{2n}{\alpha}, |\mathsf{R}_\gamma| \leq \frac{n}{\beta}\}$.

Fix $\beta=5\alpha^2$ and $\theta=5\alpha^3$, and set $r=\left\lfloor\frac n\theta\right\rfloor$, $N_1=\left\lceil\frac n\alpha\right\rceil$, $N_2=\left\lfloor\frac n\beta+3\alpha\frac n\theta\right\rfloor$. If $r=0$, then $n/\theta<1$ and every bridge in $\manyZZfewR$ belongs to $\SAB_{n,k,\theta}^{1}$, since its initial vertex is a renewal point. Thus \eqref{eq:case3_result} follows immediately. We therefore assume that $r\geq1$.
We see that for $\gamma \in \manyZZfewR$ and $Z \subset \ZigZag_\gamma$ such that $|Z| = r$, $\operatorname{Unf}_Z(\gamma) \in \SAB_{n,k,\theta}^{1}$ since each point of zig is sent to a renewal point (and $\gamma_0$ is a renewal point). 

Define the map
$$
\begin{matrix}
    \operatorname{MultiUnf} : &\manyZZfewR &\to &\mathcal{P}(\SAB_{n,k,\theta}^{1})\\
    &\gamma &\mapsto &\{\operatorname{Unf}_Z(\gamma) : Z \subset \shortZZ_\gamma , |Z| = \lfloor \frac{n}{\theta} \rfloor \},
\end{matrix}
$$
where $\shortZZ_{\gamma} := \{(i,j) \in \ZigZag_\gamma : j-i \leq \alpha\}$.

We want to use the multi-valued map principle on this map. First, we consider the size of the images. Notice that $|\ZigZag_\gamma \backslash \shortZZ_\gamma|\alpha\leq n$, and hence
$|\shortZZ_\gamma|\geq|\ZigZag_\gamma|-n/\alpha\geq n/\alpha$.
Since $|\shortZZ_\gamma|$ is an integer, we obtain
$|\shortZZ_\gamma|\geq N_1$. For each $\gamma \in \manyZZfewR$, the map $Z \to \operatorname{Unf}_Z(\gamma)$ from $\mathcal{P}(\ZigZag_\gamma) \to \SAB_n$ is one to one. Hence,
$$
|\operatorname{MultiUnf}(\gamma)| = \binom{|\shortZZ_\gamma|}{r} \geq \binom{N_1}{r}.
$$

We will now consider the size of the preimages. First, we notice that $|\mathsf{R}_\phi| \leq N_2$. Indeed, consider a single unfolding of a short zigzag $(i,j)$. Outside the range of heights
\[
y(\gamma_j)\leq y<y(\gamma_i)+2\bigl(y(\gamma_i)-y(\gamma_j)\bigr),
\]
renewal points of the unfolded bridge correspond to renewal points before the unfolding. This range contains
$3\bigl(y(\gamma_i)-y(\gamma_j)\bigr)\leq3\alpha$
heights. Thus each unfolding creates at most $3\alpha$ new renewal points. Applying this estimate successively to the $r$ unfoldings gives the claim.

As in Case 2, if $\phi \in \operatorname{MultiUnf}(\gamma)$ for a given $\phi$ and unknown $\gamma$, then it suffices to know the set of points in $\phi$ that were points of zig or zag in $\gamma$ and were unfolded by $\operatorname{MultiUnf}$. This time, we will first identify the unfolded points of zig and then identify the unfolded points of zag, knowing where the unfolding points of zig are. To find the unfolded points of zig, we only need to identify them among the renewal points of $\phi$. Hence, there are at most $\binom{|\mathsf{R}_\phi|}{r}$ possible choices. For the unfolded points of zag, we only need to know the distance to their respective points of zig. This distance is at most $\alpha$, hence there are at most $\alpha^{r}$ possible choices.
Therefore, 
$$
|\operatorname{MultiUnf}^{-1}(\phi)| \leq \binom{|\mathsf{R}_\phi|}{r} \alpha^{r} \leq \binom{N_2}{r} \alpha^{r}.
$$ 
Combining the two bounds, we get
$$
  \frac{\min_\gamma |\operatorname{MultiUnf}(\gamma)|}{\max_\phi |\operatorname{MultiUnf}^{-1}(\phi)|}
\geq \left( \frac{N_1-r}{N_2 \alpha} \right)^r,
$$
where we used that $\frac{\binom{n_1}{m}}{\binom{n_2}{m}} \geq \left(\frac{n_1 - m}{n_2}\right)^m$ for $n_1 \geq n_2 \geq m$.

By the multi-valued map principle, we therefore get
\begin{align*}
    |\SAB_{n,k,\theta}^1| 
    \geq \frac{p}{2} \frac{\min_\gamma |\operatorname{MultiUnf}(\gamma)|}{\max_\phi |\operatorname{MultiUnf}^{-1}(\phi)|} |\SAB_{n,k,\kappa}^\ell|
    \geq \frac{p}{2} \left( \frac{N_1-r}{N_2 \alpha} \right)^r |\SAB_{n,k,\kappa}^\ell|.
\end{align*}
With these choices of $\beta$ and $\theta$ (for which $\beta\leq\theta$), we have
$$
\frac{N_1-r}{N_2 \alpha}\geq \frac{5}{4} - \frac{1}{4 \alpha^2}>1,$$ 
since we assumed $\alpha \geq 8$. This gives
\begin{equation}\label{eq:case3_result}
    |\SAB_{n,k,\theta}^1| \geq \frac{p}{2} |\SAB_{n,k,\kappa}^\ell|.
\end{equation}

We will now iterate these cases to obtain the desired bound. Let $n \in \N$. Let $\theta_0=5\alpha^3$, $\alpha_0=\alpha$, and define recursively $\alpha_m$, $\beta_m$ and $\theta_m$ by
\begin{equation}\label{eq:unfolding_recursion}
    \theta_{m} 
    = 
    \begin{cases}
    \frac{1}{2}q\left(\frac{1}{10 \rho(\alpha_{m-1})}\right) & \text{if } \alpha_{m-1} \geq 8\\
    0 & \text{otherwise}
    \end{cases},
    \quad \alpha_m = \left(\frac{\theta_m}{5}\right)^{1/3} ,\quad \beta_m = 5\alpha_m^2.
\end{equation}
Notice that for all $m \geq 1$, $\theta_{m} = \bar q(\theta_{m-1})$. Let $M = M(\theta_1,\bar q)$. If $M = 0$ then the statement of the proposition is trivial, since its left-hand side is zero. Suppose therefore that $M \geq 1$.

Let us prove that the sequences $\alpha_m$ and $\theta_m$ are non-increasing in $m$. Indeed, $q(x) \leq x$ and $\rho(x) \geq 4 /x $, hence for any $m\geq 1$ such that $\theta_{m} \geq 1$, by \eqref{eq:unfolding_recursion}, we have
\begin{equation}\label{eq:theta_decreasing}
    \theta_m \leq \frac{1}{20 \rho(\alpha_{m-1})} \leq  \frac{\alpha_{m-1}}{80} = \frac{1}{80} \left(\frac{\theta_{m-1}}{5}\right)^{1/3}.
\end{equation}
By the definition of $M$, we have $\theta_{M+1}\geq3^M\geq1$, hence by a backward induction in $m$ one can prove that $\theta_m$ is non-increasing in $m$ and that for each $m \in [0;M+1]$, $\theta_m \geq 1$. One can then (\ref{eq:theta_decreasing}) to get
\begin{equation}\label{eq:alpha_decreasing}
    \alpha_m = \left(\frac{\theta_{m}}{5}\right)^{1/3} \leq \left(\frac{1}{100 \rho (\alpha_{m-1})}\right)^{1/3} \leq \frac{1}{20} \left(20 \alpha_{m-1}\right)^{1/3}
\end{equation}
for each $m \in [1;M]$. This gives in particular that $\alpha_m$ is also non-increasing. Combining \eqref{eq:theta_decreasing}, \eqref{eq:alpha_decreasing} and $\theta_{M+1} \geq 1$, one gets that for each $m \in [0;M]$, $\alpha_m \geq 80$.

We will iterate Case 2 until either Case 1 or Case 3 happens, or until we have iterated Case 2 enough times. For each $m \in [1;M]$, let $p_m \in (0,1)$. At each step $m$, we use the parameters $\alpha =\alpha_m$, $\beta =\beta_m$ and $\theta= \theta_m$ and $p = p_m$.

Fix $k = 3^M/2$. We start with the set $\SAB_{n,k,2k}^{2k}$. If Case 2 happens $m$ times, we get that
\begin{equation}\label{eq:case2_iterated}
|\SAB_{n,k,2k}^{\lfloor \frac{2k}{3^m}\rfloor}| \geq A(m) |\SAB_{n,k,2k}^{2k}|,
\end{equation}
where
$$
A(m) = \prod_{i=1}^{m} (1-p_i)e^{-\rho(\alpha_i) n}.
$$

Suppose that Case 2 happens $M$ consecutive times. Then, we can stop iterating and \eqref{eq:case2_iterated} reads as
$$
|\SAB_{n,k,2k}^{1}| 
\geq A(M) |\SAB_{n,k,2k}^{2k}|.
$$

Suppose now that this is not the case, which means that there exists $m \in [1;M]$ such that Case 1 or Case 3 happens at the $m$-th step. If Case 1 happens, then we can deduce from \eqref{eq:case1_result} and \eqref{eq:case2_iterated} that
$$
|\SAB_{n,k,\beta_{m}}^{1}| 
\geq \frac{p_{m}}{2} A(m-1) |\SAB_{n,k,2k}^{2k}|.
$$
Similarly, if Case 3 happens, then we get from \eqref{eq:case3_result} and \eqref{eq:case2_iterated} that
$$
|\SAB_{n,k,\theta_{m}}^{1}| 
\geq \frac{p_{m}}{2} A(m-1) |\SAB_{n,k,2k}^{2k}|.
$$

Therefore, regardless of which case happens, there exists $m \in [1;M+1]$ such that
\begin{equation}\label{eq: p_m}
|\SAB_{n,k,\theta_{m}}^{1}| 
\geq \frac{p_{m}}{2} A(m-1) |\SAB_{n,k,2k}^{2k}|,
\end{equation}
with $p_{M+1} = 1$ for convenience. For such an $m$, we have
$$
\P_{\SAB_n}\left(\spn(\gamma) \geq \frac{n}{k}\right) = \frac{|\SAB_{n,k}|}{|\SAB_n|} \leq \frac{|\SAB_{n,k,2k}^{2k}|}{|\SAB_n|}
$$
and
$$
\frac{|\SAB_{n,k,\theta_{m}}^{1}|}{|\SAB_n|} \leq \P_{\SAB_n}\left(|R_\gamma| \geq \frac{n}{\theta_{m}}\right),
$$
hence 
$$
\P_{\SAB_n}\left(\spn(\gamma) \geq \frac{n}{k}\right) \leq \frac{2}{p_mA(m-1)} \P_{\SAB_n}\left(|R_\gamma| \geq \frac{n}{\theta_{m}}\right).
$$

Define $P_m= p_{m} \prod_{i=1}^{m-1} (1-p_i)$. One can easily see that the values of $p_i$ can be chosen to give any desired positive values to the $P_m$ (as long as the condition $\sum_{m=1}^{M+1} P_m = 1$ is satisfied). For any $m$ satisfying \eqref{eq: p_m}, Lemma \ref{lemma:bound_renewal} implies that
$$
\P_{\SAB_n}\left(\spn(\gamma) \geq \frac{n}{k}\right)
\leq \frac{2}{P_m} \exp\left(\sum_{i=1}^{m-1} \rho(\alpha_i) n-2 \rho(\alpha_{m-1})n + C_{HW} \sqrt{n}\right).
$$
We can now choose 
$$
P_m = \frac{\exp\left(\sum_{i=1}^{m-1} \rho(\alpha_i) n-2 \rho(\alpha_{m-1})n\right)}{\sum_{m=0}^M\exp\left(\sum_{i=1}^{m} \rho(\alpha_i) n-2 \rho(\alpha_{m})n\right)}
$$ 
to make the bound independent of $m$:
$$
\P_{\SAB_n}\left(\spn(\gamma) \geq \frac{n}{k}\right)
\leq 2\sum_{m=1}^{M+1} \exp\left(\sum_{i=1}^{m-1} \rho(\alpha_i) n-2 \rho(\alpha_{m-1})n + C_{HW} \sqrt{n}\right).
$$

Using the inequality $1+\log(t)\leq 2\sqrt{t}$ we have 
\begin{equation}\label{eq:upper_bound_rho}
\rho(x)=\frac{4}{x}\left(1+\log\left(\frac{x}{4}\right)\right)\leq \frac{4}{\sqrt{x}}    
\end{equation}
By \eqref{eq:upper_bound_rho}, the bound $\rho(x)\geq 4/x$ for $x\geq 4$, and \eqref{eq:alpha_decreasing}, we have, for every $m\in[1;M]$,
\[
\frac{\rho(\alpha_{m-1})}{\rho(\alpha_m)}
\leq \frac{\alpha_m}{\sqrt{\alpha_{m-1}}}
\leq \frac{1}{20\sqrt{\alpha_m}}
\leq \frac{1}{80\sqrt{5}}
\leq \frac13,
\]
where we used that $\alpha_m\geq 80$. In particular, $\rho(\alpha_m)\geq 3\rho(\alpha_{m-1})$. We claim that, for every $m\in[0;M]$,
\[
2\rho(\alpha_m)-\sum_{i=1}^{m}\rho(\alpha_i)
\geq (m+2)\rho(\alpha_0).
\]
Indeed, this holds for $m=0$. Moreover, increasing $m$ by one increases the left-hand side by $\rho(\alpha_{m+1})-2\rho(\alpha_m)\geq \rho(\alpha_m)\geq \rho(\alpha_0)$, so the claim follows by induction.

Setting $a=\rho(\alpha_0)n>0$, the preceding probability estimate therefore gives
\begin{align*}
\P_{\SAB_n}\left(\spn(\gamma)\geq \frac{n}{k}\right)
&\leq 2e^{C_{HW}\sqrt{n}}
\sum_{m=0}^{M}
\exp\left(
\left(\sum_{i=1}^{m}\rho(\alpha_i)-2\rho(\alpha_m)\right)n
\right)\\
&\leq 2e^{C_{HW}\sqrt{n}}
\sum_{m=0}^{\infty}e^{-(m+2)a}\\
&=\frac{2e^{-2a+C_{HW}\sqrt{n}}}{1-e^{-a}}\\
&\leq \frac{e^{-a/2+C_{HW}\sqrt{n}}}{1-e^{-a/2}}\\
&=\frac{\exp\left(-\frac{\rho(\alpha_0)}{2}n+C_{HW}\sqrt{n}\right)}
{1-\exp\left(-\frac{\rho(\alpha_0)}{2}n\right)}.
\end{align*}
Here, the last inequality follows from
$1-e^{-a}=(1-e^{-a/2})(1+e^{-a/2})$ and
$2e^{-3a/2}\leq 1+e^{-a/2}$ for $a>0$.
\end{proof}

\subsection{Deducing Theorem \ref{thm: main theorem}}\label{sec:a_analysis}
In this section, we will prove Theorem \ref{thm: main theorem}, assuming that Theorem \ref{thm:stepb} holds. To do so, we first translate Proposition \ref{prop:unfolding_tree} into a result on walks and then rewrite $M(x,f)$ in a more convenient form, up to some explicit multiplicative constant. Combined with Theorem \ref{thm:stepb}, this will allow us to prove Theorem \ref{thm: main theorem}.

\begin{proposition}\label{prop:stepa_M}
    For all $n \in \N^*$ and $\alpha\geq 8$,
\[
\P_{\SAW_n}\left(\max_{k\in[1;n]}\|\gamma_k\|\geq \frac{2d^{1/2}n}{3^{M(\bar q(5\alpha^3),\bar q)}}\right)\leq2d\frac{\exp\left(-\frac{\rho(\alpha)}{2}n+2C_{HW}\sqrt{n}\right)}{1-\exp\left(-\frac{\rho(\alpha)}{2}n\right)}.
\]
\end{proposition}

The proof is the same as the one of \cite[Theorem 1.1]{DCHammond2013}. We include it nonetheless for completeness:

\begin{proof}
    Notice that $\P_{\SAW_n}$ is invariant under permutations and reflections of the Euclidean coordinates. Therefore,
    $$
    \P_{\SAW_n}\left(\max_{k\in[1;n]}\|\gamma_k\| \geq \frac{2d^{1/2} n}{3^{M(\bar q(5\alpha^3),\bar q)}}\right) \leq 2d \P_{\SAW_n}\left(\max_{k\in[1;n]}y(\gamma_k) \geq \frac{2n}{3^{M(\bar q(5\alpha^3),\bar q)}}\right).
    $$
    
    Hammersley and Welsh introduced a map from the set of self-avoiding walks to the set of self-avoiding bridges that unfolds walks in such a way that each bridge has at most $e^{C_{HW}\sqrt{n}}$ preimages. Note that $\max_{k\in[1;n]}y(\gamma_k)$ can only increase after unfolding, and $\max_{k\in[1;n]}y(\gamma_k) = y(\gamma_n)$ as soon as $\gamma$ is a bridge. Therefore,
    \begin{align*}
        \P_{\SAW_n}\left(\max_{k\in[1;n]}y(\gamma_k) \geq \frac{2n}{3^{M(\bar q(5\alpha^3),\bar q)}}\right) 
        &\leq e^{C_{HW}\sqrt{n}}\P_{\SAB_n}\left(y(\gamma_n) \geq \frac{2n}{3^{M(\bar q(5\alpha^3),\bar q)}} \right).
    \end{align*}
    Combining the two inequalities above with Proposition \ref{prop:unfolding_tree} gives the desired result.
\end{proof}

We will conclude this section by an analysis of the function $M(x,f)$. Let $f$ be a nondecreasing function for which there exists a constant $c_0>0$ such that $x-f(x) \geq c_0$ for all $x\geq 1$. Define the function $f^*$ by the following recurrence equation:

\begin{equation*}
    f^*(x) := 
    \begin{cases}
        0 &\text{if } x \leq 1,\\
        1 + f^*(f(x))  &\text{if } x > 1.
    \end{cases}
\end{equation*}

\begin{proposition}\label{prop:iter_equiv}
    Let $f$ be defined as above. Suppose that there exists $c > 1$ such that $f(x) \leq x/c$ for all $x \geq 1$. Then one has for all $x > 1$,
    \begin{equation}\label{eq:iter_equiv}
        \left\lfloor \frac{f^*(x) -1}{1 + \frac{\log 3}{\log c}} \right \rfloor \leq M(x,f) \leq f^*(x)-1.
    \end{equation}
\end{proposition}

\begin{proof}
    Let $x>1$. Denote $N = f^*(x)$ and $M = M(x,f)$. For all $m\geq 0$, define $x_m = f^{\circ m}(x)$, and note that $N=\min\{m\geq 0: x_m\leq 1\}$. By definition,
    $$
    x_N\leq 1, \quad x_{N-1}> 1.
    $$
    
    For the left inequality, let $m = \left\lfloor \frac{N -1}{1 + \frac{\log3}{\log c}} \right \rfloor$. The assumption on $f$ gives $x_{m+1} \leq x_{m}/c$. We can iterate backwards to get
    $$
    x_m \geq c^{N-1-m} x_{N-1} > c^{N-1-m}.
    $$
    Now, the choice of $m$ implies $N-1 -m\geq m \frac{\log3}{\log c}$. Therefore $x_m > 3^{m}$,
    which implies that $m \leq M$.

    For the right inequality, one can see that for any $m \geq 1$, if $x_m \geq 3^m$, then $x_m > 1$ and $m \leq N-1$. Therefore, $M \leq N-1$.
\end{proof}

Note that a better upper bound on $f$ means a sharper lower bound in \eqref{eq:iter_equiv}. In particular, this implies that if $f(x) = o(x)$ and $f(x)\to\infty$ as $x$ goes to infinity, then $M(x,f) = f^*(x)(1+o(1))$. Note that, by the same argument as \eqref{eq:theta_decreasing}, $\bar q(x)$ satisfies the hypothesis of Proposition \ref{prop:iter_equiv} with $c = 80 \cdot5^{1/3}$.
We can therefore deduce the following theorem.
\begin{theorem}\label{thm:stepa_inverse}
For all $n\in\N^*$ and $\alpha\geq 8$,
\[
\P_{\SAW_n}\left(\max_{k\in[1;n]}\|\gamma_k\| \geq \frac{6d^{1/2}n}{C^{\bar q^*(5\alpha^3)-2}}\right)\leq 2d\frac{\exp\left(-\frac{\rho(\alpha)}{2}n+2C_{HW}\sqrt{n}\right)}{1-\exp\left(-\frac{\rho(\alpha)}{2}n\right)},
\]
where
$C=\exp\left(\frac{\log3\log c}{\log3+\log c}\right)$ and $c=80\cdot5^{1/3}$.
\end{theorem}
\begin{proof}
Let $N=\bar q^*(5\alpha^3)$. If $N\leq 2$, then the event on the left-hand side is empty, since its displacement threshold is larger than $n$. Suppose therefore that $N\geq 3$. Then $\bar q(5\alpha^3)>1$ and
$\bar q^*(\bar q(5\alpha^3))=N-1$. Proposition \ref{prop:iter_equiv} gives
\[
M(\bar q(5\alpha^3),\bar q)\geq\left\lfloor
\frac{N-2}{1+\frac{\log3}{\log c}}\right\rfloor.
\]
Using that $\lfloor x\rfloor\geq x-1$ and the definition of $C$, we obtain
\[3^{M(\bar q(5\alpha^3),\bar q)}\geq3^{\frac{N-2}{1+\frac{\log3}{\log c}}-1}=\frac13 C^{N-2}.
\]
Consequently,
\[
\frac{2d^{1/2}n}{3^{M(\bar q(5\alpha^3),\bar q)}}\leq\frac{6d^{1/2}n}{C^{N-2}},
\]
and the result follows from Proposition \ref{prop:stepa_M}.
\end{proof}

We are now ready to prove Theorem \ref{thm: main theorem} assuming Theorem \ref{thm:stepb}.

\begin{proof}[Proof of Theorem \ref{thm: main theorem} assuming Theorem \ref{thm:stepb}]
We have that $q(h(x)) > x$ for all $x \geq 1$. Define the function $\bar h(x) = 5\rho^{-1}(1/(10h(2x)))^3$ and set

\[
L(v):=\left\lceil\log_C\left(\frac{6d^{1/2}}{v}\right)\right\rceil+2,
\qquad
T(v):=\left(\frac{\bar h^{\circ L(v)}(1)}{5}\right)^{1/3}.
\]
Since $\bar q(\bar h(x)) > x$ for all $x \geq 1$, applying Theorem \ref{thm:stepa_inverse} for $\alpha = T(v)$ with the bound $\rho(\alpha) \geq 4/\alpha$ gives
\begin{equation}\label{eq:bound_delta_alpha}
\P_{\SAW_n}\left(\max_{k\in[1;n]}\|\gamma_k\| \geq v n \right) \leq 2d \frac{\exp\left(- \frac{2n}{T(v)} + 2 C_{HW} \sqrt{n}\right)}{1-\exp\left(-\frac{2n}{T(v)} \right)}.   
\end{equation}
Since $\rho$ is nonincreasing on $[4,\infty)$, the inverse of any nonincreasing majorant of $\rho$ majorizes $\rho^{-1}$ and thus by \eqref{eq:upper_bound_rho} $\rho^{-1}(y)\leq (4/y)^2$ for all $y\in (0,1]$, from which we easily see that $\bar h(x)\leq h_1(x)$ for all $x\in [1,\infty)$ where $h_1(x):=5 \cdot 40^6 \cdot h(2x)^6$. Since $h_1$ is nondecreasing, we also have 
\[
\bar h^{\circ \left(n\right)}(x)\leq h_1^{\circ \left(n\right)}(x)
\]
for all $x\in [1,\infty)$ and $n\in\NN$.
Theorem \ref{thm:stepb} gives $h(x)=e^{e^{O(x)}}$, and thus also $h_1(x)=e^{e^{O(x)}}$. Since we then have $h_1(x)\leq \exp^{\circ(c)}(x)$ for some fixed integer $c>0$, $h_1^{\circ \left(n\right)}(1)$ is bounded from above by a tower of exponentials of height $O(n)$.
Since $L(v)=O(\log(1/v))$ as $v\to 0$, this implies that $T(v)$ is bounded from above by a tower of exponentials of height $O(\log(1/v))$.

For any $d > 0$, $C \geq 0$ and $T >0$, the condition $n\ge4C^2T^2+2T\log(1+2d)$ gives
$$
2d \frac{\exp\left(- \frac{n}{T} + 2 C \sqrt{n}\right)}{1-\exp\left(-\frac{2n}{T} \right)} \leq 1.
$$
Indeed, one can bound the numerator and denominator separatly. For the numerator, the AM-GM inequality on $2 C \sqrt{n}$ followed by the condition gives 
$$
-\frac{n}{T} +2C \sqrt{n} \leq -\frac{n}{2T} + 2 C^2 T \leq -\log(1+2d),
$$
which in turn gives
$$
2d \exp\left(- \frac{n}{T} + 2 C \sqrt{n}\right) \leq \frac{2d}{1+2d}.
$$
For the denominator, the condition implies $n \geq 2T \log(1+2d)$, which gives
$$
1-\exp\left(-\frac{2n}{T} \right) \geq 1-(1+2d)^{-4} \geq 1-\frac{1}{1+2d} = \frac{2d}{1+2d}.
$$

Combined with \eqref{eq:bound_delta_alpha}, this proves the theorem.
\end{proof}

\begin{remark}\label{rem:optimal_T}
    If instead we assumed the bound $h(x)=O(x^p)$, which is believed to be optimal for some value of $p\geq 1$ that depends on the dimension $d$, then the same argument as above allows us to take $h_1(x)=(C_1x)^{r}$ for $r=6p$ and some $C_1>0$, and by a simple induction we have 
	\[
	h_1^{\circ \left(n\right)}(1)=e^{\log(C_1)r \frac{r^n-1}{r-1} }=e^{O(r^n)}.
	\]
    This implies that $T(v) = \exp(O(v^{-\gamma}))$ for some explicit $\gamma>0$, and one can use that to get that with high probability the displacement is at most $O(n/(\log n)^\varepsilon)$ for some explicit $\varepsilon > 0$. On the other hand, improving the bound to $h(x)=e^{O(x)}$ leads in just the same manner to tower exponential bounds for $T(v)$.
\end{remark}

\section{Proof of Theorem \ref{thm:stepb}}
\label{sec:b}

In this section, we prove Theorem \ref{thm:stepb}, namely that $\eisab[|\gamma|]=\infty$ with the explicit rate 
\[
\eisab\left[|\gamma|^{\leq e^{e^{O(v)}}}\right]>v \text{ for $v\geq 1$}.
\] 
As highlighted in Section \ref{sec:structure}, our proof combines two quantitative ingredients to reach a contradiction, reflecting the very same structure of the qualitative argument given in \cite[Section 4]{DCHammond2013}. The first ingredient is a polynomial lower bound on the probability of the event $\{0\in C_n\}$ under the assumption that the truncated expected length of an irreducible bridge does not grow too fast, and the second one is a polynomial upper bound on the probability of the event $\{0\in C_n\}$ under a similar assumption. Since both conclusions cannot hold at the same time, this means that one of the two premises does not hold, which gives the desired result.

Before we proceed with the proofs, let us give the definition of a diamond point. We say that $k$ is a \emph{diamond point} of a bridge $\gamma$ if $y(\gamma_i)>|x(\gamma_i)|$ for every $i>k$, and $y(\gamma_i)<-|x(\gamma_i)|$ for every $i<k$. We also write $\{\gamma_{r_k}\in D_n\}$ to denote the event $y(\gamma_i)>|x(\gamma_i)|$ for every $i\in [r_k+1;r_{k+n}]$, and $y(\gamma_i)<-|x(\gamma_i)|$ for every $i\in [r_{k-n};r_k-1]$. 

\begin{figure}[ht]
    \centering
    \begin{tikzpicture}[x=0.38cm,y=0.23cm]
    \draw[cone] (0,0)--(4.5,4.5)--(1,8)--(-3.5,3.5)--cycle;
    \draw[cone] (1,8)--(4,11)--(-1,16)--(-4,13)--cycle;
    \draw[cone] (-1,16)--(-4,19)--(2,25)--(5,22)--cycle;

\draw[walk] 
plot coordinates {(0,0) (0,1) (0,2) (-0.5,2) (-0.5,3) (0,3) (0,4) (-1,4) (-1,5) (0,5) (1,5) (1,4) (1,3) (2,3) (2,4) (3,4) (3,5) (2,5) (2,6) (1,6) (1,7) (1,8) (1,9) (1,10) (0,10) (0,11) (1,11) (2,11) (2,12) (1,12) (0,12) (-1,12)  (-1,11) (-1.5,11) (-1.5,12) (-1.5,13) (-1,13) (0,13) (0,14) (-1,14) (-1,15) (-1,16) (-1,17) (-1,18) (-1.5,18) (-1.5,19) (-1.5,20) (-1,20) (-1,21) (0,21) (0,20) (0,19) (1,19) (1,20) (1,21) (2,21) (2,22) (3,22) (3,23) (2,23) (2,24) (2,25)};

\foreach \y in {(0,0), (1,8), (-1,16), (2,25)}{ \node[point] at \y {};}
\end{tikzpicture}
    \caption{An example of a bridge with four diamond points highlighted in blue.}
    \label{fig: diamond points}
\end{figure}

\subsection{Polynomial lower bound for diamond points}

The desired lower bound for the probability of a diamond point is obtained by making the argument of \cite[Proposition 4.2]{DCHammond2013} quantitative. The argument of Duminil-Copin and Hammond uses the law of large numbers to show the existence of diamond points under the assumption of finite expected length for irreducible bridges. Our argument avoids the use of the law of large numbers and is instead constructive.

This alternative strategy for establishing the existence of diamond points is based on deriving a metastable version of the law of large numbers from Etemadi's proof \cite{Etemadi1981}. While the argument presented below follows a different route, it is also inspired by Etemadi's proof. Indeed, using a Borel-Cantelli-type of argument we can bound the length of each irreducible bridge $\gamma^i$ by $ci$ for some constant $c>0$. Conditioning on this event, this allows us to work with the symmetric random variables $Z_j=\sum_{i=1}^j \tfrac{x(\gamma^i)}{i}$, whose summands are bounded in absolute value by $c$, and in order to show that $\gamma$ stays in a cone, it suffices to show that all $Z_j$ are uniformly bounded with not too small probability. This is where our argument deviates from Etemadi's proof.

The first ingredient we use to show the latter is Levy's maximal inequality for sums of independent symmetric random variables \cite[Theorem 8.3.2]{AL06}. We include a short proof for the reader's convenience.

\begin{lemma}\label{lem: maximal inequality}
Let $(Y_n)$ be a sequence of independent random variables symmetrically distributed around $0$. For every $n\geq 1$ and $x>0$ we have 
\[
\PP\left(\max_{l\in [1;n]}\left\vert Y_1+\ldots+Y_l\right\vert>x\right)\leq 2\PP\left(\left\vert Y_1+\ldots+Y_n \right\vert>x\right).
\]    
\end{lemma}
\begin{proof}
Let $\tau:=\min\{l\in [1;n]: |Y_1+\ldots+Y_l|>x  \}$. Note that conditionally on $\{\tau=k\}$, $Y_{k+1}+\ldots+Y_n$ is independent from $Y_1,\ldots, Y_k$ and symmetrically distributed around $0$. Hence, conditionally on $\{\tau=k\}$, with probability at least $1/2$, $Y_{k+1}+\ldots+Y_n$ has the same sign as $Y_1+\ldots+Y_k$. Moreover,
\[
\{\tau=k, (Y_1+\ldots+Y_k)(Y_{k+1}+\ldots+Y_n)\geq 0\} \subset \{|Y_1+\ldots+Y_n|>x\},
\]
and the events on the left are disjoint as $k$ varies.
Thus,
\begin{align*}
\PP\left(\max_{l\in [1;n]}\left\vert Y_1+\ldots+Y_l\right\vert>x \right)&=\sum_{k=1}^n \PP\left(\tau=k\right) \\ &\leq 2 \sum_{k=1}^n \PP\left(\tau=k, (Y_1+\ldots+Y_k)(Y_{k+1}+\ldots+Y_n)\geq 0\right) \\& \leq 2\PP\left(\left\vert Y_1+\ldots+Y_n \right\vert>x\right).
\end{align*}
This completes the proof.
\end{proof}

We next use Levy's maximal inequality to uniformly bound sums of independent symmetric random variables by a given constant $\varepsilon>0$, at a cost which is exponential in the variance of the sum. The main idea of the proof is the use of a reflection argument which is applied each time the sum becomes larger than $\varepsilon/2$. The variance of the sum allows us to control the number of times the sum becomes larger than $\varepsilon/2$.

\begin{lemma}\label{lem: max sum control}
Take some real $C>0$ and integer $M\geq 1$, and fix $\varepsilon>0$. Let also $Y_1,\ldots,Y_M$ be mutually independent random variables which are all symmetrically distributed around $0$ such that almost surely, $|Y_i|\leq \varepsilon/16$ for every $i=1,\ldots,M$, and
\[
\sum_{i=1}^M \mathbb{E}(Y_i^2)\leq C.
\]
Then
\[
\PP\left(\max_{1\leq j\leq M}
\bigg|\sum_{i=1}^j Y_i\bigg|<\varepsilon
\right)\geq 4^{-\lceil 64C/\varepsilon^2\rceil}.
\]
\end{lemma}

\begin{proof}
First, let us partition the interval $[1;M]$ into disjoint intervals $[a_1;b_1],\ldots,[a_N;b_N]$ such that for every $j=1,\ldots,N-1$,
\begin{equation}\label{eq: variance decomposition}
\frac{\varepsilon^2}{64}< \sum_{i=a_j}^{b_j} \EE(Y_i^2)\leq \frac{\varepsilon^2}{32},
\end{equation}
and for $j=N$,
\begin{equation}\label{eq: variance decomposition 2}
\sum_{i=a_N}^{b_N} \EE(Y_i^2)\leq \frac{\varepsilon^2}{32}.    
\end{equation}
These bounds will allow us to handle each $\max_{a_j\leq k\leq b_j}
\left|\sum_{i=a_j}^k Y_i\right|$.
We can construct these intervals greedily. Indeed, let $a_1=1$ and define $S_1:=\{1\leq k \leq M: \sum_{i=a_1}^k \EE(Y_i^2)> \varepsilon^2/64\}$. Now let $b_1:=\min S_1$ if $S_1$ is non-empty, and set $b_1:=M$ otherwise. Observe that if $S_1\neq \emptyset$, then case \eqref{eq: variance decomposition} happens for $j=1$, since, by definition, $\sum_{i=a_1}^{b_1-1} \EE(Y_i^2)\leq \varepsilon^2/64$, and $\EE(Y_{b_1}^2)\leq \varepsilon^2/64$ thanks to our assumption that $|Y_{b_1}|\leq \varepsilon/8$. On the other hand, if $S_1=\emptyset$, then case \eqref{eq: variance decomposition 2} happens. If $b_1=M$, then we stop. Otherwise, we set $a_2:=b_1+1$, and we define $b_2$ similarly to how $b_1$ was defined by working in the interval $[a_2;M]$ instead of $[a_1;M]$. We continue until we have exhausted the interval $[1;M]$. It follows that \eqref{eq: variance decomposition} and \eqref{eq: variance decomposition 2} hold.

Let us also observe that 
\begin{equation}\label{eq: N bound}
N\leq \Bigl\lceil\frac{64C}{\varepsilon^2}\Bigr\rceil,    
\end{equation}
since by \eqref{eq: variance decomposition} and our assumption,
\[
(N-1)\frac{\varepsilon^2}{64}< \sum_{i=1}^M \EE(Y_i^2)\leq C.
\]
We next handle each $\max_{a_j\leq k\leq b_j}
\bigg|\sum_{i=a_j}^k Y_i\bigg|$. By Chebyshev's inequality and \eqref{eq: variance decomposition}, \eqref{eq: variance decomposition 2},
\[
\PP\left(\bigg|\sum_{i=a_j}^{b_j} Y_i\bigg|\geq \varepsilon/2\right)\leq \frac{4}{\varepsilon^2}\sum_{i=a_j}^{b_j} \EE(Y_i^2)\leq 1/8.
\]
By Lemma~\ref{lem: maximal inequality},
\begin{equation}\label{eq: a_j b_j bound}
\PP\left(\max_{a_j\leq k\leq b_j}
\bigg|\sum_{i=a_j}^k Y_i\bigg|< \varepsilon/2\right)\geq 3/4.
\end{equation}
Now for $x,y\in (-\varepsilon/2,\varepsilon/2)$, let us condition $\sum_{i=1}^{a_j-1} Y_i=x$ and $|\sum_{i=a_j}^{b_j} Y_i|=|y|$, and note that either
$\left|x+y\right|<\varepsilon/2$ or $\left|x-y\right|<\varepsilon/2$.
Since conditionally on $\sum_{i=1}^{a_j-1} Y_i=x$, $\Big|\sum_{i=a_j}^{b_j} Y_i\Big|=|y|$, $\sum_{i=a_j}^{b_j} Y_i$ is a symmetric random variable,
\[
\PP\left(\bigg|\sum_{i=1}^{b_j} Y_i\bigg|<\varepsilon/2 \ \  \Big| \ \  \sum_{i=1}^{a_j-1} Y_i=x, \bigg|\sum_{i=a_j}^{b_j} Y_i\bigg|=|y|\right)\geq 1/2.
\]
By averaging and using independence,
\[
\PP\left(\bigg|\sum_{i=1}^{b_j} Y_i\bigg|< \varepsilon/2 \ \ \Big| \ \ \mathcal{E}_{j-1}, \bigg|\sum_{i=a_j}^{b_j} Y_i\bigg|<\varepsilon/2 \right)\geq 1/2,
\]
where
\[
\mathcal{E}_{j-1}:=\left\{\forall \, m=1\ldots,j-1\max_{a_m\leq k\leq b_m}
\bigg|\sum_{i=a_m}^k Y_i\bigg|< \varepsilon/2, \bigg|\sum_{i=1}^{b_m} Y_i\bigg|< \varepsilon/2\right\}.
\]
Using a union bound, independence, and \eqref{eq: a_j b_j bound},
\begin{align*}
\PP\left(\max_{a_j\leq k\leq b_j}
\bigg|\sum_{i=a_j}^k Y_i\bigg|< \varepsilon/2, 
\bigg|\sum_{i=1}^{b_j} Y_i\bigg|< \varepsilon/2 \ \ \Big| \ \ \mathcal{E}_{j-1} \right)&\geq 3/4-\PP\left(
\bigg|\sum_{i=a_j}^{b_j} Y_i\bigg|< \varepsilon/2, 
\bigg|\sum_{i=1}^{b_j} Y_i\bigg|\geq \varepsilon/2 \ \ \Big| \ \ \mathcal{E}_{j-1} \right)\\ 
&\geq 3/4-\PP\left( \bigg|\sum_{i=1}^{b_j} Y_i\bigg|\geq \varepsilon/2 \ \ \Big| \ \ \mathcal{E}_{j-1}, \bigg|\sum_{i=a_j}^{b_j} Y_i\bigg|< \varepsilon/2 \right)\\
&\geq 1/4.
\end{align*}
Thus, iterating the latter bound and using \eqref{eq: N bound},
\[
\PP\left(\mathcal{E}_N \right)\geq 4^{-\lceil 64C/\varepsilon^2\rceil}.
\]
Finally, on the event $\mathcal{E}_N$, for every $1 \leq k \leq M$, letting $m$ be the largest integer such that $a_m\leq k$,
\[
\bigg|\sum_{i=1}^k Y_i\bigg|\leq \bigg|\sum_{i=1}^{b_{m-1}} Y_i\bigg|+\bigg|\sum_{i=a_m}^{k} Y_i\bigg|< \varepsilon.
\]
This completes the proof.
\end{proof}

In the above lemma, we chose to bound the probability that the maximum of the sum remains smaller than $\varepsilon$ in terms of the variance of the sum, while it is also possible to bound this probability in terms of the expectation of the sum of $|Y_i|$. The reason why we prefer to use the variance here will become clear in Lemma~\ref{lem: bounds} below, as the variance is the natural quantity that arises when we apply Lemma~\ref{lem: max sum control} to bridges.

Next, we show how to adapt Lemma~\ref{lem: max sum control} to averages of random variables.

\begin{lemma}\label{lem: averages}
Fix $M\geq 1$ and $\varepsilon>0$, and suppose that $X_1,\ldots,X_M$ are mutually independent random variables which are symmetrically distributed around $0$ such that $|X_i|\leq i\varepsilon/32$ almost surely for every $i=1,\ldots,M$, and
\[
\sum_{i=1}^M\frac{\EE[X_i^2]}{i^2}\leq C.
\]
Then writing $S_i:=X_1+\ldots+X_i$ we have
\[
\PP\left(\forall i=1,\ldots,M\, \frac{|S_i|}{i}\leq \varepsilon
\right)\geq 4^{-\lceil 256C/\varepsilon^2\rceil}.
\]
\end{lemma}

\begin{proof}
Let $Z_0=0$ and for every $j\geq1$ define
\[
Z_j:=\sum_{i=1}^j \frac{X_i}{i}.
\]   
Note that the random variables $X_1,X_2/2,\ldots,X_M/M$ satisfy the assumptions of Lemma~\ref{lem: max sum control} on parameter $\varepsilon/2$, hence we get
\[
\PP( \forall i=1,\ldots,M\, |Z_i|<\varepsilon/2)\geq 4^{-\lceil 256C/\varepsilon^2\rceil}.
\]   
To relate the latter bound to the desired bound of the lemma, we use summation by parts. More precisely, note that $X_i=i(Z_i-Z_{i-1})$. Then,
\[
\sum_{i=1}^j X_i=\sum_{i=1}^j i Z_i-\sum_{i=0}^{j-1} (i+1) Z_i=jZ_j-\sum_{i=0}^{j-1}Z_i.
\]
By the triangle inequality, it follows that on the event $\{\forall i=1,\ldots,M\, |Z_i|<\varepsilon/2\}$,
\[
|S_j|\leq j|Z_j|+\sum_{i=0}^{j-1}|Z_i|< \varepsilon j.
\]
This completes the proof.
\end{proof}

In what follows, we now fix two sequences $(X_n)$ and $(Y_n)$ of random variables. Furthermore, we assume that these are related as follows:
\begin{itemize}
	\item the pairs $(X_n,Y_n)_{n\geq 1}$ are i.i.d.,
	\item $(X_1,Y_1)$ and $(-X_1,Y_1)$ have the same distribution,
	\item $|X_1|\leq Y_1$ a.s.,
	\item there exists $b\geq 0$ with $p:=\PP(Y_1\leq b)>0$ but $\PP(Y_1\leq b,X_1>0)=0$.
\end{itemize} 
For $M\geq 1$ and $\varepsilon>0$ let
\[
A_{M,\varepsilon}:=\{\forall\, i=1,\ldots,N_\varepsilon-1\, Y_i\leq b, \forall i = N_\varepsilon,\ldots,M\, Y_i\leq i\varepsilon/32\}
\]
for $N_{\varepsilon}:=\max\{\lceil 32b/\varepsilon\rceil,1\}$. We next show that conditionally on $A_{M,\varepsilon}$ and under the assumption that the truncated expectation of $Y_1$ is bounded by a given constant, $X_n/n$ satisfies the assumption of Lemma~\ref{lem: max sum control}.

\begin{lemma}\label{lem: bounds}
Fix $M\geq 1$, $\varepsilon\in (0,1)$ and $v>0$.
Then whenever $\EE[Y_1^{\leq M}]\leq v$ we have
\[
\sum_{n=1}^M \frac{\EE[X_n^2\mid A_{M,\varepsilon}]}{n^2}\leq 2v/p.
\]
\end{lemma}

\begin{proof}
Let $N:=N_{\varepsilon}$ and write $c_n:=n\varepsilon/32$. Note that conditionally on $A_{M,\varepsilon}$, $X_n=0$ for $n=1,2,\ldots,N-1$ since $\PP(Y_n\leq b,X_n>0)=0$, so we can restrict the sum to $n\in [N;M]$. Following a broadly similar argument to the proof of Lemma 5.20 in \cite{Klenke2020}, we first observe that for $n\in [N;M]$, $|X_n|\leq Y_n$ a.s.\ and mutual independence of the $(Y_n)_{n\geq 1}$ gives us
\[
\PP(|X_n|>t\mid A_{M,\varepsilon})\leq \PP(Y_n>t\mid A_{M,\varepsilon})=\PP(Y_n>t\mid Y_n\leq c_n)=\frac{\PP(t<Y_1\leq c_n)}{\PP(Y_1\leq c_n)}\leq \frac{\PP(Y_1^{\leq M}>t)}{p},
\]
where the last inequality follows from the fact that $c_n\geq b$ and thus $\PP(Y_1\leq c_n)\geq \PP(Y_1\leq b)=p$, and also $c_n\leq n\leq M$. Therefore
\begin{align*}
\EE[X_n^2\mid A_{M,\varepsilon}]&=\int_0^\infty 2t\, \PP(|X_n|>t\mid A_{M,\varepsilon})\, dt\\
&\leq 2\int_0^\infty \mathbf{1}_{t<c_n}\,t\, \PP(Y_n>t\mid Y_n\leq c_n)\, \mathrm{d}t \leq \frac{2}{p} \int_0^\infty \mathbf{1}_{t<c_n}\,t\, \PP(Y_1^{\leq M}>t)\, \mathrm{d}t.
\end{align*}
It therefore follows that
\begin{align*}
\sum_{n=1}^M \frac{\EE[X_n^2\mid A_{M,\varepsilon}]}{n^2}&\leq 2p^{-1}\sum_{n=N}^M n^{-2}\int_0^\infty \mathbf{1}_{t<c_n}\,t\, \PP(Y_1^{\leq M}>t)\, \mathrm{d}t.
\end{align*}
Hence, because the integrand is nonnegative, we can change the order of summation so the right-hand side becomes
\[
2p^{-1}\int_0^\infty \left(\sum_{n=N}^M n^{-2}\mathbf{1}_{t<c_n}\right)\,t\, \PP(Y_1^{\leq M}>t)\, dt
\] 
and since
\[
\sum_{n=N}^M n^{-2}\mathbf{1}_{t<c_n}\leq \sum_{n=1}^\infty n^{-2}\mathbf{1}_{t<c_n}\leq \sum_{n>2t}n^{-2}\leq \frac{1}{t}
\]
the integral is bounded by
\[
2p^{-1}\int_0^\infty  \PP(Y_1^{\leq M}>t)\, dt=2p^{-1}\EE[Y_1^{\leq M}]\leq 2v/p
\]
and this completes the proof.
\end{proof}

We now assume that we have an explicit tail rate function $\kappa:(0,1]\rightarrow \NN$ for $Y_1$ i.e.\ satisfying $\PP(Y_1 \geq\kappa(\lambda))\leq \lambda$ for all $\varepsilon\in (0,1)$. This allows us to construct a lower bound on the probability of $A_{M,\varepsilon}$.

\begin{lemma}\label{lem:Abound}
Fix $\varepsilon\in (0,1)$ and $M\geq N:=N_{\varepsilon}$. Then whenever $\EE[Y_1^{\leq \kappa(1/M)}]\leq v$ we have
\[
\PP\left(A_{M,\varepsilon}\right)\geq p^{N-1}\exp\left(-32vp^{-1}\varepsilon^{-1}-p^{-1}\right).
\]
\end{lemma}

\begin{proof}
We will use that for $x\in [0,1)$ we have $\log(1-x)\geq -x/(1-x)$, which follows from applying the standard inequality $t\geq 1-e^{-t}$ for $t=\log(1-x)$. Writing $t_n:=n\varepsilon/32$ as before, and noting again that $\PP(Y_1\leq t_n)\geq p$ for $n\in [N,M]$, by independence we have
\begin{align*}
\PP\left(A_{M,\varepsilon}\right)&=p^{N-1} \exp\left(\sum_{n=N}^M\log\left(1-\PP\left(Y_n> t_n\right)\right)\right)\\
&\geq p^{N-1}\exp\left(\sum_{n=N}^M-\frac{\PP\left(Y_1> t_n\right)}{\PP\left(Y_n\leq t_n\right)}\right)\geq p^{N-1}\exp\left(-p^{-1}\sum_{n=N}^M\PP\left(Y_1> t_n\right)\right).
\end{align*}
Now, since
\begin{align*}
\PP\left(Y_1> t_n\right)&\leq\PP\left(\kappa(M^{-1})\geq Y_1>t_n\right)+\PP\left(Y_1>\kappa(M^{-1})\right) \\ &\leq \PP\left(Y_1^{\leq \kappa(M^{-1})}>t_n\right)+\frac{1}{M},
\end{align*}
we can conclude
\[
\sum_{n=N}^M\PP\left(Y_1> t_n\right)\leq \sum_{n=1}^\infty \PP\left(Y_1^{\leq \kappa(1/M)}>t_n\right)+\sum_{n=1}^M \frac{1}{M}\leq \frac{32}{\varepsilon}\EE[Y_1^{\leq \kappa(1/M)}]+1\leq \frac{32v}{\varepsilon}+1
\]
and the result follows.
\end{proof}

We now collect everything together into the following key technical lemma that will be required to establish our lower bound for diamond points.

\begin{lemma}\label{lem:technical}
Let $(X_n,Y_n)_{n\geq 1}$ be i.i.d.\ with $(X_1,Y_1)$ and $(-X_1,Y_1)$ identically distributed and $|X_1|\leq Y_1$ a.s. Suppose that $b\geq 0$ is such that $p:=\PP(Y_1\leq b)>0$ and $\PP(Y_1\leq b,X_1>0)=0$. Finally, let $\kappa$ be a tail rate for $Y_1$. Then for any $\varepsilon\in (0,1)$ and $M\geq N_\varepsilon:=\max\{\lceil 32b/\varepsilon\rceil,1\}$, whenever $\EE[Y_1^{\leq\max\{M,\kappa(1/M)\}}]\leq v$ we have
\begin{align*}
&\PP\left(\forall i\in [1;N_\varepsilon)\, Y_i\leq b, \forall i\in [N_\varepsilon,M]\left(|X_1+\ldots+X_i|\leq i\varepsilon,Y_i\leq i\varepsilon/32\right)\right)\\
&\geq 4^{-\lceil 512vp^{-1}\varepsilon^{-2}\rceil}p^{N_\varepsilon-1}\exp\left(-32vp^{-1}\varepsilon^{-1}-p^{-1}\right).
\end{align*}
\end{lemma}

\begin{proof}
Fixing $\varepsilon\in (0,1)$ and $M\geq N_\varepsilon$, we apply Lemma \ref{lem: averages} to $X_1,\ldots,X_M$ under the conditional measure $\PP(\ \cdot \mid A_{M,\varepsilon})$, where the assumptions of that lemma are satisfied conditionally on $A_{M,\varepsilon}$ due to the independence of $(X_n,Y_n)_{n\geq 1}$, $(X_n,Y_n)\overset{d}{=}(-X_n,Y_n)$ and $|X_n|\leq Y_n$ a.s., together with Lemma \ref{lem: bounds}. From this we obtain
\[
\PP(\forall i\in [1;M]\, |X_1+\ldots+X_i|\leq i\varepsilon\mid A_{M,\varepsilon})\geq 4^{-\lceil 512vp^{-1}\varepsilon^{-2}\rceil}
\]
and the result then follows by Lemma \ref{lem:Abound}.
\end{proof}

Though we do not use it, we note that the following quite natural result on i.i.d.\ and symmetrically distributed random variables corresponds to the special case of our technical lemma with $Y_n:=|X_n|$.

\begin{corollary}
Let $(X_n)_{n\geq 1}$ be i.i.d.\ and symmetrically distributed around $0$. Suppose that $p:=\PP(X_1=0)>0$ and let $\kappa$ be a tail rate for $|X_1|$. Then for any $\varepsilon\in (0,1)$ and $M\geq 1$, whenever $\EE[|X_1|^{\leq \max\{M,\kappa(1/M)\}}]\leq v$ we have
\begin{align*}
&\PP\left(\max_{1\leq i\leq M}\, \frac{|X_1+\ldots+X_i|}{i}\leq \varepsilon\right)\geq 4^{-\lceil 512vp^{-1}\varepsilon^{-2}\rceil}\exp\left(-32vp^{-1}\varepsilon^{-1}-p^{-1}\right).
\end{align*}
\end{corollary}

We now arrive at our main result of this section, namely our quantitative variant of Proposition 4.2 of \cite{DCHammond2013}.

\begin{theorem}\label{thm:contradiction-a}
If $\eisab[|\gamma|]<\infty$ then $\pisabn(0\in C)>0$. Moreover, we can explicitly define a constant $K_1\geq 1$ such that for all $v\geq 1$ and $n\geq 64$, whenever 
\[
\eisab\left[|\gamma|^{\leq \phi(1/n)}\right]\leq v
\]
then $\pisabn\left(0\in C_n\right)\geq e^{-K_1v}$, for $\phi$ as in \eqref{eq: phi def}.
\end{theorem}

\begin{proof}
We apply Lemma \ref{lem:technical} on the measure $\pisabn$ for $X_n:=x(\gamma_{r_n})-x(\gamma_{r_{n-1}})$, $Y_n:=|\gamma^n|$, and $\kappa(\lambda):=\phi(\lambda)$ for $\phi$ as defined in \eqref{eq: phi def}, noting that $\pisabn(|\gamma^1|\leq 1)=\pisab(|\gamma|\leq 1)=\pisab(|\gamma|=1)=\mu_c^{-1}$ and $\pisabn(|\gamma^1|\leq 1,x(\gamma_{r_1})>0)=\pisab(|\gamma|\leq 1,x(\gamma)>0)=0$, and so we set $b:=1$ and $p:=\mu_c^{-1}$. In the case $\varepsilon=1/2$, we get that for any $n\geq 64$, whenever $\eisabn[|\gamma_1|^{\leq \phi(1/n)}]=\eisab[|\gamma|^{\leq\phi(1/n)}]\leq v$ (also observing that $\phi(\lambda)\geq 1/\lambda$ and thus $\max\{M,\phi(1/M)\}=\phi(1/M)$) then
\begin{align*}
&\pisabn\left(\forall i\in [1;64)\, |\gamma^i|=1, \forall i\in [64,n](|x(\gamma_{r_i})|\leq i/2, |\gamma^i|\leq i/64) \right)\\
&\geq 4^{-\lceil 2^{11}v\mu_c\rceil}\mu_c^{-63}\exp\left(-64v\mu_c-\mu_c\right)\geq e^{-K_1v}
\end{align*}
for suitable $K_1\geq 1$ dependent on $\mu_c$. Let $A$ be the event in the probability measure. We use the abbreviations $x_i,y_i:=x(\gamma_{r_i}),y(\gamma_{r_i})$. For any $i$ we have $y_i+x_i\geq i-|x_i|$ using the fact that the height of any irreducible bridge is at least one and thus $y_i\geq i$. Thus for $k\in [r_{i},r_{i+1}]$ it follows that 
\[
y(\gamma_k)+x(\gamma_k)\geq i-|x_i|-|\gamma^{i+1}|
\]
accounting for the worst case in which the path $\gamma_{[r_i,k]}$ is east-west. By an analogous argument we have $y(\gamma_k)-x(\gamma_k)\geq i-|x_i|-|\gamma^{i+1}|$. Therefore on $A$, noting that $\gamma^i$ must be a north-south edge for each $i\in [1,64]$, if $k\in [64,r_n]$ then $k\in [r_{i-1},r_{i}]$ for some $i\in [65,n]$, and so we have
\[
\min((y+x)(\gamma_k),(y-x)(\gamma_k))\geq (i-1)-|x_{i-1}|-|\gamma^{i}|> 0,
\]
where for the last inequality we use that $|x_{i-1}|+|\gamma^{i}|\leq (i-1)/2+i/64< i-1$. On the other hand, for $k\leq 64$ we have $\gamma_k=(0,k,0,\ldots,0)$. We have thus shown that $A\subseteq \{0\in C_{n}\}$, and the second assertion follows. The first part is then immediate, since $\eisab[|\gamma|]<\infty$ implies that there exists some $v\geq 1$ such that $\eisab[|\gamma|^{\leq M}]\leq v$ for all natural numbers $M$, and so in particular
\[
\pisabn(0\in C)=\inf_{n\in\NN}\pisabn(0\in C_n)\geq e^{-K_1v}>0.\qedhere
\]
\end{proof}

\subsection{A modified stickbreak argument}

We now provide a modification of the proof of Theorem 2.5 of \cite{DCHammond2013} (cf.\ pages 22--25 of that paper) that combines with our computational Theorem \ref{thm:contradiction-a}. Once again our argument avoids the use of the law of large numbers that was used to handle certain unlikely events, which we perform instead by quantitative arguments. We follow a similar multi-valued map construction and the main quantitative gain comes from a tighter control on the number of pre-images using a geometric construction.

\begin{figure}[ht]
\centering
\begin{tikzpicture}[x=0.25cm,y=0.20cm]

\foreach \y in {0,8,16,24}{\draw[cone] (0,\y)--(4,{\y+4})--(0,{\y+8})--(-4,{\y+4})--cycle;}
\draw[walk] 
plot coordinates {(0,0) (0,1) (0,2) (1,2) (1,3) (2,3) (2,4) (2,5) (1,5) (0,5) (-1,5) (-1,6) (0,6) (0,7) (0,8)(0,9) (0,10) (-1,10) (-1,11) (-1,12) (-2,12) (-2,13) (-1,13) (0,13) (0,12) (1,12) (1,11) (2,11) (2,12) (2,13) (1,13) (1,14) (0,14) (0,15) (0,16) (0,17) (0,18) (-0.5,18) (-0.5,19) (0,19) (0,20) (-1,20) (-1,19) (-2,19) (-2,20) (-2,21) (-1,21) (0,21) (1,21) (1,22) (0,22) (0,23) (0,24) (0,25) (0,26) (-1,26) (-1,27) (-1,28) (0,28) (0,27) (1,27) (2,27) (2,28) (2,29) (1,29) (0,29) (0,30) (0,31) (0,32)};

\foreach \y in {0,16,32}{\node[point] at (0,\y) {};}
\node[chosen] at (0,8) {};
\node[chosen] at (0,24) {};

\begin{scope}[xshift=2cm]

\draw[cone] (12,0)--(16,4)--(12,8)--(8,4)--cycle;
\draw[cone] (12,8)--(16,4)--(20,8)--(16,12)--cycle;
\draw[cone] (20,8)--(24,4)--(28,8)--(24,12)--cycle;
\draw[cone] (28,8)--(32,12)--(28,16)--(24,12)--cycle;

\draw[walk] plot coordinates {(12,0) (12,1) (12,2) (13,2) (13,3) (14,3) (14,4) (14,5) (13,5) (12,5) (11,5) (11,6) (12,6) (12,7) (12,8) (13,8) (14,8) (14,9) (15,9) (16,9) (16,10) (17,10) (17,9) (17,8) (16,8) (16,7) (15,7) (15,6) (16,6) (17,6) (17,7) (18,7) (18,8) (19,8) (20,8) (21,8) (22,8) (22,8.5) (23,8.5) (23,8) (24,8) (24,9) (23,9) (23,10) (24,10) (25,10) (25,9) (25,8) (25,7) (26,7) (26,8) (27,8) (28,8) (28,9) (28,10) (27,10) (27,11) (27,12) (28,12) (28,11) (29,11) (30,11) (30,12) (30,13) (29,13) (28,13) (28,14) (28,15) (28,16)};

\foreach \p in 
{(12,0),(20,8),(28,16)}{\node[point] at \p {};}

\node[chosen] at (12,8) {};
\node[chosen] at (28,8) {};
\end{scope}  
\end{tikzpicture}
\caption{The stickbreak operation at two selected diamond points highlighted in red. The middle section is rotated clockwise, while the prefix and suffix
are left unchanged up to translation.}
\label{fig: stickbreak operation}
\end{figure}

\begin{lemma}
\label{lem:stickbreak}
For any $n\geq n_0:=6^{10}$ we have
\[
\eisab\left[|\gamma|^{\leq \phi\left(n^{-11/10}\right)}\right]\leq n^{1/10}\implies \pisabn(0\in C_n) <7n^{-1/20}
\]
for $\phi$ as in \eqref{eq: phi def}.
\end{lemma}
\begin{proof}
Fix $n\geq n_0$, and let $M:=\phi\left(n^{-11/10}\right)$. Then
\[
\pisabn\left(\exists i\in [1;n]\left(|\gamma^i|>M\right)\right)\leq n\cdot \pisab\left(|\gamma|>M\right)\leq n\cdot n^{-11/10}=n^{-1/10}
\]
and thus
\[
\pisabn\left(\forall i\in [1;n]\left(|\gamma^i|\leq M\right)\right)\geq 1-n^{-1/10}.
\]

Let us set $m:=n/3+n^{4/5}-1$. To simplify the notation we assume that $n/3$ and $m$ are integers, as otherwise one can take integer parts. Let $\mathcal E_n$ denote the event that among the renewal points $r_i$ with $i\in \{ n/3,n/3+1,\ldots,m\}$, for
at least $n^{7/10}$ of them $\gamma_{r_i}\in D_n$. We first estimate $\pisabn(\mathcal{E}_n)$, which we will later use to estimate $\pisabn(0\in C_n)$.

Define the events
\[
\mathcal A_n := \Bigl\{\max_{r_{n/3}<k\leq r_m}|x(\gamma_k)-x(\gamma_{r_{n/3}})|<\tfrac n4 \Bigr\},
\qquad
\mathcal B_n := \Bigl\{ \sum_{i=1}^n|\gamma^i|\le n^{6/5} \Bigr\}.
\]
By the assumed moment bound,
\[
\eisabn\!\left[\mathbbm{1}_{\forall\, i\in [1;n]\, |\gamma^i|\leq M}
\max_{r_{n/3}<k\leq r_m}|x(\gamma_k)-x(\gamma_{r_{n/3}})| \right] \le n^{9/10}.
\]
Hence, by a union bound and Markov's inequality,
\[
\pisabn(\mathcal A_n)
\ge 1-5n^{-1/10}.
\]
Similarly, $\eisabn [ \sum_{i=1}^n|\gamma^i|^{\leq M}]\leq n^{11/10}$, which together with another application of Markov's inequality and a union bound yields
\[
\pisabn(\mathcal B_n)
\ge 1-2n^{-1/10}.
\]
Let $\mathcal E_n' := \mathcal E_n \cap \mathcal A_n \cap \mathcal B_n$.
By the union bound,
\begin{equation}\label{eq:union bound}
\pisabn(\mathcal E_n) \le \pisabn(\mathcal E_n') + 7n^{-1/10}.
\end{equation}

Our aim now is to bound $\pisabn(\mathcal E_n')$. To do this, we will use the multi-valued map principle. We define a map $F$ as follows.
Let $\gamma\in\mathcal E_n'$ and choose two distinct diamond points with renewal indices
$i<j$, both lying between $n/3$ and $m$. Write the walk as a concatenation $
\gamma = \gamma^{(1)} \circ \gamma^{(2)} \circ \gamma^{(3)}$,
where $\gamma^{(1)}$ is the walk up to the renewal point $r_i$, $\gamma^{(2)}$ is the segment
from $r_i$ to $r_j$, and $\gamma^{(3)}$ is the walk from $r_j$ to $r_n$. Define a
new walk $F(\gamma,i,j)$ by rotating $\gamma^{(2)}$ by $\pi/2$ and gluing it back at $\gamma_{r_i}$, and gluing $\gamma^{(3)}$ at the $r_j$ index of this new walk, namely
\[
F(\gamma,i,j):= \gamma^{(1)} \circ R(\gamma^{(2)}) \circ \gamma^{(3)},
\]
where $R(\gamma^{(2)})$ denotes the rotation of $ \gamma^{(2)}$ by $\pi/2$ around its starting point. See Figure~\ref{fig: stickbreak operation}.

By our assumption that $\gamma_{r_i}$ and $\gamma_{r_j}$ are diamond points, $F(\gamma,i,j)$ is a self-avoiding walk. Note that $y(\gamma_{r_{n/3}})\geq n/3$, because each irreducible component has height at least $1$, and $\gamma\in \mathcal{A}_n$, hence we have that $y(F(\gamma,i,j))_k>0$ for all $k\geq 1$. Similarly,
\begin{equation*}
\begin{aligned}
y(F(\gamma,i,j))_{r_n}&\geq \max_{k\leq r_{n/3}}y(F(\gamma,i,j))_k-\max_{r_{n/3}<k\leq r_m}|x(\gamma_k)-x(\gamma_{r_{n/3}})|+n-m\\&>\max_{k\leq r_{n/3}}y(F(\gamma,i,j))_k+\max_{r_{n/3}<k\leq r_m}|x(\gamma_k)-x(\gamma_{r_{n/3}})|\geq \max_{k\leq r_{m}}y(F(\gamma,i,j))_k,
\end{aligned}
\end{equation*}
where here we emphasise that we denote by $r_i$ the $i$th renewal point of $\gamma$, not of $F(\gamma,i,j)$.
Since this is also clearly true for $r_m<k<r_n$, it follows that $F(\gamma,i,j)$ is a bridge.

Let us write $\mathcal{S}$ to denote the set of all $F(\gamma,i,j)$, where $\gamma\in \mathcal{E}'_n$ and $i,j$ are as above. We will bound $\pisabn(\mathcal{E}'_n)$ in terms of $\pisabn(\mathcal{S})$. To this end, first note that for each $\gamma\in \mathcal{E}'_n$, both $\gamma$ and $F(\gamma,i,j)$ have the same length, hence $\pisabn(\gamma)=\pisabn(F(\gamma,i,j))$. Moreover, on the one hand, each $\gamma\in \mathcal{E}'_n$ has at least $
\binom{n^{7/10}}{2}$ possible images. On the other hand, we claim that each $\gamma'\in \mathcal{S}$ has at most $n^{6/5}$ preimages.

Indeed, by the definition of the event
$\mathcal B_n$, there are at most $n^{6/5}$ possible choices for the step $\ell$ where the second stickbreak was applied.
Fix such a choice of $\ell$, and rotate by $\pi/2$ the part of $\gamma'$ from $\gamma'_{\ell}$ to $\gamma'_{r_n}$,
thereby reconstructing an intermediate self-avoiding walk $\tilde\gamma$. In $\tilde\gamma$, due to the first stickbreak operation, there must be a step $k$ such that 
\begin{enumerate}
    \item $(y+x)(\tilde\gamma)_r<(y+x)(\tilde\gamma)_k$ and $(y-x)(\tilde\gamma)_r<(y-x)(\tilde\gamma)_k$ for every $r<k$,
    \item $(y+x)(\tilde\gamma)_r>(y+x)(\tilde\gamma)_k$ and $(y-x)(\tilde\gamma)_r<(y-x)(\tilde\gamma)_k$ for every $r>k$.
\end{enumerate} 
It is not hard to see that $k$ is the only point of $\tilde\gamma$ with these properties due to the strict inequalities in the definitions. Indeed, for $r<k$,  we have $(y+x)(\tilde\gamma)_r<(y+x)(\tilde\gamma)_k$, which means that $r$ does not satisfy property $(2)$. For $r>k$, we have $(y+x)(\tilde\gamma)_r>(y+x)(\tilde\gamma)_k$, which means that $r$ does not satisfy property $(1)$. This uniquely determines the first step where the stickbreak operation happened as the unique point that satisfies properties $(1)$ and $(2)$, and we can, in turn, uniquely reconstruct the original $\gamma$ from $\tilde\gamma$ by applying a $-\pi/2$ rotation to the subwalk of $\tilde\gamma$ from $k$ to $r_n$. This proves the claim. See Figure~\ref{fig: unique stickbreak point}.

\begin{figure}[ht]
\centering
\begin{tikzpicture}[x=0.58cm,y=0.58cm]

\draw[cone] (-5,-5)--(5,5);
\draw[cone] (-5,5)--(5,-5);
  
\draw[cone] (0,-3)--(2.5,-0.5);
\draw[cone] (0,-3)--(2.5,-5.5);

\draw[past] 
plot coordinates {(-1,-6) (-1,-5) (1,-5) (1,-4) (0,-4) (0,-3) (-1,-3) (-1,-2) (0,-2) (0,-1) (0,0)};

\draw[future] 
plot coordinates { (0,0) (1,0) (2,0) (2,1) (3,1) (3,-1) (4,-1) (4,1) (5,1) (5,0) (6,0)
};

\node[point] at (0,-3) {};
\node[anchor=east,font=\small] at (-0.05,-3.35) {$r$};
\node[chosen] at (0,0) {};
\node[anchor=south east,font=\small] at (0.4,0.2) {$k$};
\end{tikzpicture}
\caption{An example illustrating the uniqueness of $k$. The future of the walk lies in the cone to the right of $k$, while its past lies in the cone below $k$. If $r$ also had this property, then the future of the walk from $r$ would have to lie in the cone to its right. This is impossible because $k$, which belongs to the future of $r$, lies outside that cone.}
\label{fig: unique stickbreak point}
\end{figure}

Combining the above, we obtain that 
\[
\pisabn(\mathcal{E}'_n)\leq \frac{1}{\binom{n^{7/10}}{2}}\sum_{\gamma\in \mathcal{E}'_n}\sum_{i,j} \pisabn(F(\gamma,i,j))\leq \frac{n^{6/5}}{\binom{n^{7/10}}{2}} \pisabn(\mathcal{S})\leq \frac{n^{6/5}}{\binom{n^{7/10}}{2}},
\]
where we used that we have at least $\binom{n^{7/10}}{2}$ choices for $i,j$ in the first inequality, and that each $F(\gamma,i,j)\in \mathcal{S}$ appears at most $n^{6/5}$ times in the sum to get the second inequality.
Using \eqref{eq:union bound}
we conclude that
\[
\pisabn(\mathcal E_n)\leq 11 n^{-1/10}.
\]
Now note that by reflection symmetry and linearity of expectation,
\[
\eisabn(|n/3 \leq i\leq m: \gamma_{r_i}\in D_n|)\geq n^{4/5} \pisabn(0\in C_n)^2,
\]
so that 
\[
\eisabn(n^{4/5}-|n/3 \leq i\leq m: \gamma_{r_i}\in D_n|)\leq n^{4/5} (1-\pisabn(0\in C_n)^2).
\]
By applying Markov's inequality to the non-negative random variable $n^{4/5}-|n/3 \leq i\leq m: \gamma_{r_i}\in D_n|$, we obtain 
\[
\pisabn(\mathcal{E}^c_n)\leq \frac{n^{4/5} (1-\pisabn(0\in C_n)^2)}{n^{4/5}-n^{7/10}},
\]
hence
\[\pisabn(0\in C_n)^2\leq 1-\frac{n^{4/5}-n^{7/10}}{n^{4/5}}\pisabn(\mathcal{E}^c_n)\leq 1-(1-n^{-1/10})(1-11n^{-1/10})< 12n^{-1/10}.
\]
The desired bound follows.
\end{proof}

The following result is an immediate consequence of the modified stickbreak argument, which in combination with Theorem \ref{thm:contradiction-a} gives us our main result.

\begin{theorem}\label{thm:contradiction-b}
If $\eisab[|\gamma|]<\infty$ then $\pisabn(0\in C)=0$. Moreover, for any $v>0$ and $\delta\in (0,1]$, whenever 
\[
\eisab\left[|\gamma|^{\leq \phi(N_{v,\delta}^{-11/10})}\right]\leq v,
\]
then $\pisabn\left(0\in C_{N_{v,\delta}}\right)<\delta$, where $N_{v,\delta}:=\max\{(7/\delta)^{20},v^{10}\}$ and $\phi$ is as in Lemma \ref{lem:markov}.
\end{theorem}

\begin{proof}
The second part follows directly from Lemma \ref{lem:stickbreak}, setting $n:=N_{v,\delta}$ (note that $n\geq 7^{20}\geq 6^{10}$ since $\delta\leq 1$). For the first part, simply note that fixing $\delta>0$ we obtain $\pisabn(0\in C)\leq\pisabn\left(0\in C_{N_{v,\delta}}\right)<\delta$, and the claim follows since $\delta$ was arbitrary.
\end{proof}

We are now ready to prove Theorem \ref{thm:stepb}.

\begin{proof}[Proof of Theorem \ref{thm:stepb}]
Fix $v\geq 1$ and define $N:=N_{v,\delta}$ as in Theorem \ref{thm:contradiction-b} for $\delta:=e^{-K_1v}$ where $K_1\geq 1$ is as in Theorem \ref{thm:contradiction-a}. It is easy to see that $N=(7e^{K_1v})^{20}$ and thus
\[
\phi(N^{-1})\leq \phi(N^{-11/10})\leq e^{CN^{11/10}}=e^{C(7e^{K_1v})^{22}}\leq e^{e^{K_2v}}
\]
where $C$ is as in Lemma \ref{lem:markov} and $K_2\geq 1$ is some suitable constant independent of $v$. Now assume that
\[
\eisab\left[|\gamma|^{\leq e^{e^{K_2v}}}\right]\leq v.
\]
By Theorem \ref{thm:contradiction-a} applied to $n:=N\geq 64$ we have $\pisab(0\in C_N)\geq e^{-K_1v}=\delta$, and by Theorem \ref{thm:contradiction-b} we have $\pisab(0\in C_N)<\delta$. Therefore the assumption is false, and since $v\geq 1$ was arbitrary the result follows.
\end{proof}

We end this section by proving Theorem~\ref{thm: b_n}, which is a reformulation of Theorem~\ref{thm:stepb}.

\begin{proof}[Proof of Theorem~\ref{thm: b_n}]
We will prove that 
\begin{equation}\label{lemma1_1}
\irrE[|\gamma|^{\leq T}] \sum_{k = 0}^{T-1} b_k \mu_c^{-k} \leq 2T.    
\end{equation}
Then the desired result follows from Theorem~\ref{thm:stepb}.

By Lemma~\ref{isab_renewal}, for any $k \in \N$, $b_k \mu_c^{-k} = \infP[\gamma_k \in \mathsf{R}_\gamma]$. Summing over $k$ we obtain
\begin{equation}\label{eq:1}
\sum_{k = 0}^{T-1} b_k \mu_c^{-k} = \sum_{k = 0}^{T-1} \infP[\gamma_k \in \mathsf{R}_\gamma] = \infE[|\mathsf{R}_\gamma[0,T-1]|].
\end{equation}
Notice that 
\begin{equation*} \label{eq:2}
|\mathsf{R}_\gamma[0,T - 1]| 
\leq R:=\inf\left\{\ell \in \N : \sum_{i = 1}^\ell |\gamma^{i}|^{\leq T} \geq T \right\}.
\end{equation*}
Therefore, we can use Wald's identity to get that
\begin{equation}\label{eq:3}
\infE\left[\sum_{i = 1}^{R} |\gamma^{i}|^{\leq T} \right] = \irrE[|\gamma|^{\leq T}] \infE[R]\geq \irrE[|\gamma|^{\leq T}] \infE[|\mathsf{R}_\gamma[0,T-1]|].
\end{equation}
By definition,
$$
\sum_{i = 1}^{R - 1} |\gamma^{i}|^{\leq T} \leq T,
$$
hence
\begin{equation}\label{eq:4}
\infE\left[\sum_{i = 1}^{R} |\gamma^{i}|^{\leq T}\right] \leq 2T,
\end{equation}
where we used that the last step is at most of size $T$. Then, \eqref{lemma1_1} follows from combining \eqref{eq:1}, \eqref{eq:3} and \eqref{eq:4}. This completes the proof.   
\end{proof}

\noindent\textbf{Acknowledgements.} We would like to thank Hugo Duminil-Copin for useful discussions and for proposing this project to EA and CP in 2022, when CP was in Geneva. EA was supported by the Swiss National Science
Foundation and the NCCR SwissMAP. CP was supported by an EPSRC New Investigator Award (UKRI1019). TP was partially supported by an EPSRC New Investigator Award (EP/W035847/1).\\

\noindent\textbf{AI declaration.} Generative AI (specifically OpenAI's GPT-6) was used for proofreading and to help check some proofs for correctness. All proofs and ideas are original, and none were produced by AI. No part of the article was written by an LLM.

\bibliographystyle{plain}
\bibliography{ref}

\begin{thebibliography}{10}

\bibitem{AL06}
K.~B. Athreya and S.~N. Lahiri.
\newblock {\em Measure Theory and Probability Theory}.
\newblock Springer, New York, 2006.

\bibitem{AvigadFeferman1998}
J.~Avigad and S.~Feferman.
\newblock {{G}\"{o}del's Functional (``{D}ialectica'') Interpretation}.
\newblock In S.R. Buss, editor, {\em Handbook of Proof Theory}, volume 137 of
  {\em Studies in Logic and the Foundations of Mathematics}, pages 337--405.
  North Holland, Amsterdam, 1998.

\bibitem{LecturesSAW}
R.~Bauerschmidt, H.~Duminil-Copin, J.~Goodman, and G.~Slade.
\newblock Lectures on self-avoiding walks.
\newblock {\em Probability and Statistical Physics in Two and More Dimensions
  (D. Ellwood, CM Newman, V. Sidoravicius, and W. Werner, eds.), Clay
  Mathematics Institute Proceedings}, 15:395--476, 2012.

\bibitem{BBDCG14}
N.~R. Beaton, M.~Bousquet-M{\'e}lou, J.~de~Gier, H.~Duminil-Copin, and A.~J.
  Guttmann.
\newblock The critical fugacity for surface adsorption of self-avoiding walks
  on the honeycomb lattice is $1+\sqrt{2}$.
\newblock {\em Communications in Mathematical Physics}, 326(3):727--754, 2014.

\bibitem{BLM2013}
S.~Boucheron, G.~Lugosi, and P.~Massart.
\newblock {\em Concentration Inequalities: A Nonasymptotic Theory of
  Independence}.
\newblock Oxford University Press, Oxford, 2013.

\bibitem{PolyIns}
H.~Duminil-Copin, S.~Ganguly, A.~Hammond, and I.~Manolescu.
\newblock Bounding the number of self-avoiding walks: {Hammersley--Welsh} with
  polygon insertion.
\newblock {\em The Annals of Probability}, 48(4):1644--1692, 2020.

\bibitem{SAWEnds}
H.~Duminil-Copin, A.~Glazman, A.~Hammond, and I.~Manolescu.
\newblock On the probability that self-avoiding walk ends at a given point.
\newblock {\em The Annals of Probability}, 44(2):955--983, 2016.

\bibitem{DCHammond2013}
H.~Duminil-Copin and A.~Hammond.
\newblock {Self-Avoiding Walk is Sub-Ballistic}.
\newblock {\em Communications in Mathematical Physics}, 324:401--423, 2013.

\bibitem{DCKY14}
H.~Duminil-Copin, G.~Kozma, and A.~Yadin.
\newblock Supercritical self-avoiding walks are space-filling.
\newblock {\em Annales de l'Institut Henri Poincar{\'e}, Probabilit{\'e}s et
  Statistiques}, 50(2):315--326, 2014.

\bibitem{DC12}
H.~Duminil-Copin and S.~Smirnov.
\newblock The connective constant of the honeycomb lattice equals $\sqrt
  {2+\sqrt 2}$.
\newblock {\em Annals of Mathematics}, 175(3):1653--1665, 2012.

\bibitem{Etemadi1981}
N.~Etemadi.
\newblock {An elementary proof of the strong law of large numbers}.
\newblock {\em Zeitschrift f\"ur Wahrscheinlichkeitstheorie und verwandte
  Gebiete}, 55(1):119--122, 1981.

\bibitem{Flory53}
P.~J. Flory.
\newblock {\em Principles of polymer chemistry}.
\newblock Cornell University Press, 1953.

\bibitem{GM20}
A.~Glazman and I.~Manolescu.
\newblock Self-avoiding walk on $\mathbb{Z}^2$ with {Yang--Baxter} weights:
  universality of critical fugacity and $2$-point function.
\newblock {\em Annales de l'Institut Henri Poincar{\'e}, Probabilit{\'e}s et
  Statistiques}, 56(4):2281--2300, 2020.

\bibitem{Goedel1958}
K.~G\"odel.
\newblock {{\"Uber eine bisher noch nicht ben\"utzte Erweiterung des finiten
  Standpunktes}}.
\newblock {\em Dialectica}, 12(3--4):280--287, 1958.

\bibitem{Hammersley1961}
J.~M. Hammersley.
\newblock On the rate of convergence to the connective constant of the
  hypercubical lattice.
\newblock {\em The Quarterly Journal of Mathematics}, 12(1):250--256, 1961.

\bibitem{Poor}
J.~M. Hammersley and K.~W. Morton.
\newblock {Poor man's Monte Carlo}.
\newblock {\em Journal of the Royal Statistical Society: Series B
  Methodological)}, 16(1):23--38, 1954.

\bibitem{HammersleyWelsh}
J.~M. Hammersley and D.~J.~A. Welsh.
\newblock Further results on the rate of convergence to the connective constant
  of the hypercubical lattice.
\newblock {\em Quarterly Journal of Mathematics}, 13(1):108--110, 1962.

\bibitem{Hammond18}
A.~Hammond.
\newblock An upper bound on the number of self-avoiding polygons via joining.
\newblock {\em The Annals of Probability}, 46(1):175--206, 2018.

\bibitem{HaSla1992}
T.~Hara and G.~Slade.
\newblock The lace expansion for self-avoiding walk in five or more dimensions.
\newblock {\em Reviews in Mathematical Physics}, 4(02):235--327, 1992.

\bibitem{HaSlaCri}
T.~Hara and G.~Slade.
\newblock {Self-avoiding walk in five or more dimensions I. The critical
  behaviour}.
\newblock {\em Communications in Mathematical Physics}, 147(1):101--136, 1992.

\bibitem{Hutchcroft2018}
T.~Hutchcroft.
\newblock The {Hammersley--Welsh} bound for self-avoiding walk revisited.
\newblock {\em Electronic Communications in Probability}, 23:Paper No. 5, 8,
  2018.

\bibitem{KestenI}
H.~Kesten.
\newblock {On the Number of Self-Avoiding Walks}.
\newblock {\em Journal of Mathematical Physics}, 4(7):960--969, 1963.

\bibitem{KestenII}
H.~Kesten.
\newblock {On the Number of Self-Avoiding Walks. II}.
\newblock {\em Journal of Mathematical Physics}, 5(8):1128--1137, 1964.

\bibitem{Klenke2020}
A.~Klenke.
\newblock {\em {Probability Theory}}.
\newblock Springer, Berlin Heidelberg, 3rd edition, 2020.

\bibitem{Kohlenbach1996}
U.~Kohlenbach.
\newblock {Analysing proofs in analysis}.
\newblock In {\em Logic: from Foundations to Applications. European Logic
  Colloquium (Keele, 1993)}, pages 225--260. Oxford University Press, 1996.

\bibitem{Kohlenbach2008}
U.~Kohlenbach.
\newblock {\em {Applied Proof Theory: Proof Interpretations and their Use in
  Mathematics}}.
\newblock Springer Monographs in Mathematics. Springer-Verlag Berlin
  Heidelberg, 2008.

\bibitem{Kohlenbach2019}
U.~Kohlenbach.
\newblock {Proof-theoretic Methods in Nonlinear Analysis}.
\newblock In B.~Sirakov, P.~Ney de~Souza, and M.~Viana, editors, {\em
  {Proceedings of ICM 2018}}, volume~2, pages 61--82. World Scientific,
  Singapure, 2019.

\bibitem{KohlenbachOliva2003}
U.~Kohlenbach and P.~Oliva.
\newblock {Proof Mining: A Systematic Way of Analysing Proofs in Mathematics}.
\newblock {\em Proceedings of the Steklov Institute of Mathematics},
  242:136--164, 2003.

\bibitem{KP26}
D.~Krachun and C.~Panagiotis.
\newblock Quantitative sub-ballisticity of self-avoiding walk on the hexagonal
  lattice.
\newblock {\em The Annals of Probability}, 54(3):1109--1125, 2026.

\bibitem{Madras14}
N.~Madras.
\newblock A lower bound for the end-to-end distance of the self-avoiding walk.
\newblock {\em Canadian Mathematical Bulletin}, 57(1):113--118, 2014.

\bibitem{MaSaSAWBook}
N.~Madras and G.~Slade.
\newblock {\em The self-avoiding walk}.
\newblock Birkh{\"a}user, Boston, 1993.

\bibitem{Neri2026}
M.~Neri.
\newblock {Oscillations in the Strong Law of Large Numbers}.
\newblock {\em Proceedings of the American Mathematical Society}, 2026.
\newblock To appear.

\bibitem{NeriOlivaPischke2026}
M.~Neri, P.~Oliva, and N.~Pischke.
\newblock {A systematic way of analysing proofs in probability theory}, 2026.
\newblock Preprint, available at \url{https://arxiv.org/abs/2604.08078}.

\bibitem{NeriPischke2024}
M.~Neri and N.~Pischke.
\newblock {Proof mining and probability theory}.
\newblock {\em Forum of Mathematics, Sigma}, 13, 2025.
\newblock Article no.\ e187, 47pp.

\bibitem{NeriPowell2025a}
M.~Neri and T.~Powell.
\newblock {On quantitative convergence for stochastic processes: Crossings,
  fluctuations and martingales}.
\newblock {\em Transactions of the American Mathematical Society, Series B},
  12:974--1019, 2025.

\end{thebibliography}

\end{document}